\title{Explicit isoperimetric constants and phase transitions in the
random-cluster model}
\author{Olle H\"{a}ggstr\"{o}m\thanks{Mathematical Statistics, G\"oteborg
University, 412 96 G\"oteborg, Sweden,
\texttt{olleh@math.chalmers.se}, 
\texttt{http://www.math.chalmers.se/\~{ }olleh/}}\
\thanks{Research partially supported by the Swedish Natural Science 
Research Council.}
\and Johan Jonasson\thanks{Department of Mathematics, Chalmers 
University of Technology, 412 96 G\"oteborg, Sweden,
\texttt{jonasson@math.chalmers.se}, 
\texttt{http://www.math.chalmers.se/\~{ }jonasson/}}
\and Russell Lyons\thanks{Department of Mathematics,  
Indiana University, Bloomington, IN 47405-5701, USA,
\texttt{rdlyons@indiana.edu}, 
\texttt{http://php.indiana.edu/\~{ }rdlyons/} and 
School of Mathematics, Georgia Institute of Technology, Atlanta, GA
30332-0160, USA,
\texttt{rdlyons@math.gatech.edu}, 
\texttt{http://www.math.gatech.edu/\~{ }rdlyons/}
}\ 
\thanks{Research partially supported by NSF grant DMS-9802663 and by
Microsoft Corp.}}

\documentclass[11pt,twoside]{article}
\usepackage{amsmath,amssymb,amsthm,amsfonts,latexsym}
\newtheorem {thm}{Theorem}[section]
\newtheorem {prop}[thm]{Proposition}

\def\Cox{\hfill \Box}
\def\Z{{\mathbb Z}}
\def\R{{\bf R}}
\def\P{{\bf P}}
\def\T{{\bf T}}
\def\Pt{{\sf Pt}}
\def\FPt{{\sf FPt}}
\def\WPt{{\sf WPt}}
\def\RC{{\sf RC}}
\def\FRC{{\sf FRC}}
\def\WRC{{\sf WRC}}
\def\C{{\cal C}}
\def\leqd{\stackrel{\cal D}{\preccurlyeq}}
\def\geqd{\stackrel{\cal D}{\succcurlyeq}}
\def\Ccon{\stackrel{\infty}{\longleftrightarrow}}

\def\Aut{{\rm Aut}}
\def\f{{\rm free}}
\def\w{{\rm wired}}

\def\ed#1{[#1]}  
\def\kk(#1){\|#1\|}  
\def\kkk(#1){\|#1\|^*}  
\def\I#1{{\bf 1}_{\{#1\}}}  
\def\st{\, ; \;}  
\def\B{{\mathbf B}}
\def\p{{\mathbf p}}
\def\pc{p_{\rm c}}
\def\pu{p_{\rm u}}
\def\pcw{\pc^\w}
\def\pcf{\pc^\f}
\def\puw{\pu^\w}
\def\puf{\pu^\f}
\def\peq{p_{{\rm F}={\rm W}}}
\def\bd{\partial}
\def\iso{\iota}
\def\bde{\bd_E}

\def\isoe{\iso_E}

\def\d#1{#1^\dagger}  

\def\fK{K^f}

\def\Xf#1{{^\f X_{#1,q}^G}}

\def\CC{{\cal C}}
\def\betac{\beta_{\rm c}}
\def\br{{\rm br}}

\def\thmenv#1#2#3{\begin{#1} \label{#1:#2} #3 \end{#1}}

\def\procl#1.#2 #3\endprocl{%
       \ifx#1t\thmenv{thm}{#2}{#3}\fi
       \ifx#1l\thmenv{lem}{#2}{#3}\fi
       \ifx#1p\thmenv{prop}{#2}{#3}\fi
       \ifx#1c\thmenv{cor}{#2}{#3}\fi
       \ifx#1d\thmenv{defn}{#2}{#3}\fi
       \ifx#1g\thmenv{conj}{#2}{#3}\fi
       \ifx#1q\thmenv{quest}{#2}{#3}\fi
       \ifx#1r\thmenv{rmk}{#2}{{\rm #3}}\fi
    }%

\def\rref#1.#2/{%
      \ifx #1sSection~\ref{sect:#2}\fi
      \ifx #1tTheorem~\ref{thm:#2}\fi  
      \ifx #1lLemma~\ref{lem:#2}\fi 
      \ifx #1cCorollary~\ref{cor:#2}\fi 
      \ifx #1pProposition~\ref{prop:#2}\fi 
      \ifx #1dDefinition~\ref{defn:#2}\fi
      \ifx #1gConjecture~\ref{conj:#2}\fi 
      \ifx #1qQuestion~\ref{quest:#2}\fi 
      \ifx #1e(\ref{eq:#2})\fi
      \ifx #1b\cite{#2}\fi
        }

\def\rlabel #1 #2{\begin{equation} \label{eq:#1} #2 \end{equation}}

\def\proof{\medbreak\noindent{\it Proof.\enspace}}

\def\proofof #1.#2 {\medbreak\noindent
     {\it Proof of \rref #1.#2/.}\enspace}

\def\Qed{$\Cox$\medbreak}

\begin{document}

\maketitle


\begin{abstract}
The random-cluster model is a dependent percolation model that has 
applications in the study of Ising and Potts models. In this paper, several
new results are obtained for the random-cluster model on nonamenable graphs
with cluster parameter $q\geq 1$. 
Among these, the main ones are the absence of percolation for the
free random-cluster measure at the critical value, and examples of planar
regular graphs with regular dual where $\pc^\f (q) > \pu^\w (q)$
for $q$ large enough.
The latter follows from
considerations of
isoperimetric constants, and we give the first
nontrivial explicit calculations of such constants. 
Such considerations are also used to prove 
non-robust phase transition for the Potts model on nonamenable
regular graphs.  

\end{abstract}


\section{Introduction} \label{sect:intro}

One of the most important and much-studied dependent percolation models 
today is the random-cluster model. It was introduced in 1972 by
Fortuin and Kasteleyn \cite{FK}, and after a decade and a half of relative 
silence, the model was revived in the late 1980s with the influential
papers by Swendsen and Wang \cite{SW}, Edwards and Sokal \cite{ES}, 
and Aizenman, Chayes, Chayes, and Newman \cite{ACCN}. Since then, the
random-cluster model
has served as a major tool in studying Ising and Potts models, and has
also been studied in its own right by several authors. This paper is
an investigation of various 
aspects of the random-cluster model on so-called nonamenable graphs, 
and we shall obtain
results that are intrinsic to the model itself, as well as some
applications to the Potts model. 
Our main results, in order of importance, are the following:
\begin{itemize}
\item We consider 
{\bf four} (possibly different) {\bf critical values} for the 
random-cluster model on a given Cayley graph, and we sort out how these
can relate to each other (Sections \ref{sect:critical} and
\ref{sect:example}).
In particular, we show that on all nonamenable planar regular graphs
with regular dual one has for $q$ large enough that
$\pc^\f (q) > \pu^\w (q)$.
For this purpose, we give the first explicit nontrivial calculations of
positive {\bf isoperimetric constants}.
\item For unimodular transitive nonamenable graphs, we show {\bf lack of
percolation\/} for the free random-cluster measure at the lower critical
value on all nonamenable quasi-transitive unimodular
graphs (\rref t.death/), thereby extending a result of
Benjamini, Lyons, Peres and Schramm \cite{BLPS, BLPScrit} that was for i.i.d.\
percolation.
\item The random-cluster model is exploited to show that the Potts model on
all nonamenable regular graphs exhibits, for entire intervals of
temperatures, phase transition but not so-called 
{\bf robust phase transition} (Section \ref{sect:robust}).
\end{itemize}

We also present a number of other results about the random-cluster measures
and about their relationships to Potts models.

We shall begin by recalling in Section \ref{sect:background} some basics
concerning random-cluster and Potts models.  
After that, we come in Section \ref{sect:critical} to the lack of 
percolation at criticality.
Proceeding to our results on the inequalities between the four critical 
values $\pc^\w$, $\pc^\f$, $\pu^\w$ and $\pu^\f$, we shall find that
most of them
are fairly immediate.
An exception is the result that 
$\pc^\f(q) > \pu^\w(q)$ can occur for some graphs and some $q$, which requires
careful analysis (carried out in Section \ref{sect:example}) of the 
random-cluster model on regular tilings of the hyperbolic plane. Parts
of this analysis are based on the Peierls-type comparison
methods of Jonasson and Steif \cite{JS} 
and Jonasson \cite{J}, and these methods are also exploited in 
Section \ref{sect:robust} to obtain our result on non-robust phase
transition 
for the Potts model. 

Since isoperimetric constants are of independent interest,
we state here our result concerning them.
A regular euclidean polygon
of $\d d$ sides has interior angles $\pi(1-2/\d d)$.
In order for such polygons to form a tessellation of the plane with $d$
polygons meeting at each vertex, we must have $\pi(1-2/\d d) = 2\pi/d$,
i.e., $1/d + 1/\d d = 1/2$, or, equivalently, $(d - 2)(\d d - 2) = 4$. 
In all three such cases, tessellations have been well known since antiquity.
In the hyperbolic plane, the interior angles can take any value in
$\big(0,\pi(1-2/\d d)\big)$, whence a tessellation exists only if  $1/d + 1/\d
d < 1/2$, or, equivalently, $(d - 2)(\d d - 2) > 4$; again, this condition
is also sufficient for the existence of a hyperbolic tesselltion,
as has been known since the 19th century.
(We remark that the cases $(d - 2)(\d d - 2) < 4$ correspond to the
spherical tessellations that arise from the five regular solids.)

Let $G=(V,E)$ be the planar graph formed by the edges and vertices of one
of these regular tessellations by polygons with $\d d$ sides and with
degree $d$ at each vertex. 
Given a finite set $K \subset V$, write $\bde K$ for the set of
edges with exactly one endpoint in $K$. We prove in \rref t.isocalc/ that
$$
\inf \left\{ {|\bde K| \over |K|} \st K \subset V \hbox{
finite and nonempty} \right\}
= (d-2) \sqrt{1 - {4 \over (d - 2)(\d d - 2)}}\,.
$$
This should be compared to the regular tree of degree $d$, where the
left-hand side is equal to $d-2$.
This formula also makes easy the task of deciding whether such a graph
satisfies one of the conditions of high nonamenability that appear in
\cite{Schon:mult}; such a condition has numerous implications for various
models on the graph.
Series and Sina{\u\i} \rref b.SeriesSinai/ were the first to consider the
Ising model on such graphs.

\section{Background} \label{sect:background}

In the following subsections,
we recall some preliminaries on random-cluster and Potts
models and on stochastic domination
and various classes of infinite graph
structures. General references for Sections 
\ref{sect:background_finite}--\ref{sect:background_infinite_Potts} 
are H\"aggstr\"om \cite{H} and
Georgii, H\"aggstr\"om and Maes \cite{GHM}, whereas for Section 
\ref{sect:background_amenability_etc}, we refer to Benjamini et al.\ 
\cite{BLPS}.  

\subsection{Random-cluster and Potts models on finite graphs}
\label{sect:background_finite}

Let $G=(V,E) = \bigl((V(G), E(G)\bigr)$ be a finite graph. An edge $e\in E$
connecting two vertices
$x,y\in V$ is also denoted $\ed{ x,y }$. 
An element $\xi$ of $\{0,1\}^E$ will
be identified with the subgraph of $G$ that has vertex set $V$ and edge
set $\{e\in E \st \xi(e)=1\}$. An edge $e$ with $\xi(e)=1$
($\xi(e)=0$) is said to be open (closed). 
Of central importance to us will be the
number of connected components of $\xi$, which will be denoted $\kk(\xi)$.
We emphasize that in the definition of $\kk(\xi)$, 
isolated vertices in $\xi$ also count as connected components. 

The {\bf random-cluster measure} $\RC:=\RC^G_{p,q}$ (sub- and superscripts
will be dropped whenever possible) with parameters $p \in [0,1]$ and $q>0$, 
is the probability measure on $\{0,1\}^E$ that to each
$\xi\in\{0,1\}^E$ assigns probability
\begin{equation} \label{eq:RC_def}
\RC(\xi):= 
\frac{q^{\kk(\xi)}}{Z} \prod_{e\in E} p^{\xi(e)}(1-p)^{1-\xi(e)} \, ,
\end{equation}
where $Z:=Z^G_{p,q}:= \sum_{\xi\in\{0,1\}^E} q^{\kk(\xi)}
\prod_{e\in E} p^{\xi(e)}(1-p)^{1-\xi(e)}$ is a normalizing constant. It
is easy to see that if $X$ is a $\{0,1\}^E$-valued random variable with
distribution $\RC$, then we have, for each $e = \ed{x,y} \in E$ and each 
$\xi\in\{0,1\}^{E\setminus\{e\}}$, that
\begin{equation} \label{eq:single_edge_cond_prob}
\RC\bigl(X(e)=1 \bigm| X(E\setminus\{e\})= \xi\bigr) =
\left\{
\begin{array}{ll}
p & \mbox{if } x\leftrightarrow y, \\
\frac{p}{p+(1-p)q} & \mbox{otherwise,}
\end{array} \right.
\end{equation}
where $x\leftrightarrow y$ is the event that there is an open path 
(i.e., a path of open edges) from $x$ to $y$ in $X(E\setminus\{e\})$.
Here, $X(E')$ denotes the restriction of $X$ to $E'$ for $E' \subseteq E$.

When $q=1$, we see that all edges are independently open and closed with
respective probabilities $p$ and $1-p$, so that we get the usual 
i.i.d.\ bond percolation model on $G$, which we refer to as {\bf
Bernoulli}($p$) percolation. All other choices of $q$ yield
dependence between the edges. Throughout the paper, {\bf we shall assume
that $q\geq 1$}. This conforms with most other studies of the random-cluster
model, and there are two reasons for doing this. First, when $q\geq 1$, 
the conditional probability in (\ref{eq:single_edge_cond_prob})
becomes increasing not only in $p$ but also in $\xi$, and this allows
a set of very powerful stochastic domination arguments that are not
available for $q< 1$. Second, it is only random-cluster measures with
$q\in\{2,3,\ldots\}$ that have proved to be useful to the analysis of 
Potts models, which we now describe. 

Given the finite graph $G$ and an integer $q\geq 2$, the $q$-state 
Potts model provides a model for picking an element 
$\omega\in \{1, \ldots, q\}^V$ in a random but correlated way.
The values $1, \ldots,q$ attainable at each vertex $x\in V$ are
called {\bf spins}. 
Fix the so-called {\bf inverse temperature} parameter $\beta\geq 0$,
and define the {\bf Gibbs measure for the $q$-state Potts model on $G$
at inverse temperature $\beta$}, denoted $\Pt:=\Pt^G_{q,\beta}$, 
as the probability measure that to each $\omega\in\{1, \ldots, q\}^V$
assigns probability 
\[
\Pt(\omega) :=
\frac{1}{Z} \exp\left( -2\beta \sum_{\ed{ x,y } \in E}
\I{\omega(x) \neq \omega(y)} \right) \, , 
\]
where $Z$ is another normalizing constant (different from the one in 
(\ref{eq:RC_def})). The main link between random-cluster and Potts
models is the following well-known result. 
\begin{prop} \label{prop:from_RC_to_Potts}
Fix the finite graph $G$, an integer $q\geq 2$, and $p\in [0,1]$. Pick
a random edge configuration $X\in \{0,1\}^E$ according to the random-cluster
measure $\RC^G_{p,q}$. Then, for each connected component $\C$ of $X$, 
pick a spin uniformly from $\{1,\ldots, q\}$, and assign this spin to 
all vertices of ${\cal C}$. Do this independently for different
connected components. The $\{1, \ldots, q\}^V$-valued random spin
configuration arising from this procedure is then distributed according
to the Gibbs measure $\Pt^G_{q,\beta}$ for the $q$-state Potts model on $G$
at inverse temperature $\beta:= -\frac{1}{2}\log(1-p)$.  
\end{prop}
This provides a way of reformulating problems about pairwise dependencies
in the Potts model into problems about connectivity probabilities in the
random-cluster model. Aizenman et al.\ \cite{ACCN} exploited such ideas
to obtain results about the phase transition behavior of the Potts model.
See \cite{H} for a list of references of other applications of the
random-cluster model to the Potts model.

\subsection{Stochastic domination and weak convergence}
\label{sect:stoch_dom}

Let $E$ be any finite or countably infinite set. (The reason for denoting it
by $E$ is that in our applications, it will be an edge set.) For two
configurations $\xi, \xi' \in \{0,1\}^E$, we write $\xi \preccurlyeq \xi'$
if
$\xi(e) \leq \xi'(e)$ for all $e\in E$. A function 
$f:\{0,1\}^E \rightarrow \R$ is said to be {\bf increasing} if $f(\xi) \leq
f(\eta)$
whenever $\xi\preccurlyeq \eta$. For two probability measures $\mu$ and
$\mu'$
on $\{0,1\}^E$, we say that $\mu$ is {\bf stochastically dominated} by $\mu'$,
writing $\mu \leqd \mu'$, if
\begin{equation} \label{eq:def_stoch_dom}
\int_{\{0,1\}^E} fd\mu \leq \int_{\{0,1\}^E} fd\mu' 
\end{equation}
for all bounded increasing $f$. By Strassen's Theorem,
this is equivalent to the existence of a coupling $\P$ of two
$\{0,1\}^E$-valued
random variables $X$ and $X'$, with respective distributions $\mu$ and $\mu'$,
such that $\P(X\preccurlyeq X')=1$. 

A useful tool for establishing stochastic domination is Holley's inequality
(see \cite{H} or \cite{GHM}). Since the conditional distribution in 
(\ref{eq:single_edge_cond_prob}) is increasing both in $\xi$ and in $p$
(recall that we consider random-cluster measures only with $q\geq 1$),
Holley's inequality applies to show that, for any finite graph $G=(V,E)$,
\[
\RC^G_{p_1,q} \leqd \RC^G_{p_2, q}
\]
whenever $p_1 \leq p_2$. Similarly we get, for conditional distributions, 
that
\begin{equation} \label{eq:cond_prob_domination}
\RC^G_{p,q} \bigl(X \in \cdot \bigm| X(E') = \xi\bigr) 
\leqd
\RC^G_{p,q} \bigl(X \in \cdot \bigm| X(E') = \xi'\bigr)
\end{equation}
whenever $E' \subseteq E$ and $\xi \preccurlyeq \xi'$. 

We shall also be considering weak convergence of probability measures
on $\{0,1\}^E$. For such probability measures $\mu_1, \mu_2, \ldots$ 
and $\mu$, we say that $\mu$ is the (weak) limit of $\mu_i$ as 
$i \rightarrow \infty$ if
$\lim_{i\rightarrow\infty}\mu_i(A) = \mu(A)$ for all cylinder events $A$. 

\subsection{The random-cluster model on infinite graphs}
\label{sect:background_infinite}

Let $G=(V,E)$ be an infinite, locally finite graph. The definition
(\ref{eq:RC_def}) of random-cluster measures
does not work in this case, because there are uncountably
many different configurations $\xi \in \{0,1\}^E$.
Instead, there are two 
other approaches to defining random-cluster measures
on infinite graphs: one via limiting procedures, and one via local 
specifications (Dobrushin-Lanford-Ruelle, or DLR, equations). 
We shall sketch the first approach, and
then explain how it relates to the second. 

Let $V_1, V_2, \ldots$ be a sequence of finite vertex sets increasing to
$V$ in the sense that $V_1 \subset V_2 \subset \ldots$ and
$\bigcup_{i=1}^\infty V_i = V$. For any finite $K \subseteq V$, define
\[
E(K) := \bigl\{ \ed{ x,y } \in E \st \, x,y \in K\bigr\},
\]
set $E_i := E(V_i)$
and note that $E_1, E_2, \ldots $ increases to $E$ in the same sense that
$V_1, V_2, \ldots$ increases to $V$. Set $\partial V_i$ to be the (inner)
boundary of $V_i$, i.e., 
\[
\partial V_i := \bigl\{v \in V_i \st \exists \, \ed{ x,y } \in E\setminus
E_i
\mbox{ with } x=v \bigr\} \, .
\] 
Also set $G_i:= (V_i, E_i)$, and let
$\FRC_{p,q}^{G, i}$ be the probability measure on $\{0,1\}^E$ corresponding
to
picking $X\in \{0,1\}^E$ by letting $X(E_i)$ have distribution
$\RC_{p,q}^{G_i}$ and setting $X(e):=0$ for all $e\in E\setminus E_i$.
Since the projection of $\FRC_{p,q}^{G, i}$ on $\{0,1\}^{E\setminus E_i}$
is deterministic, we can also view $\FRC_{p,q}^{G, i}$
as a measure on $\{0,1\}^{E_i}$, in which case it coincides with
$\RC_{p,q}^{G_i}$. By applying (\ref{eq:cond_prob_domination}) to the graph
$G_i$ with $E':= E_i \setminus E_{i-1}$ and $\xi\equiv 0$, we get that
\[
\FRC_{p,q}^{G, i-1} \leqd \FRC_{p,q}^{G, i} \, ,
\]
so that
\begin{equation} \label{eq:monotone_limit}
\FRC_{p,q}^{G, 1} \leqd \FRC_{p,q}^{G, 2} \leqd \cdots \,.
\end{equation}
This implies the existence of a limiting probability measure $\FRC^G_{p,q}$
on $\{0,1\}^E$. 
This limit is independent of the choice of $\{V_i\}_{i=1}^\infty$, and
we call it the random-cluster measure on $G$ with {\bf free boundary 
condition} (hence the ${\sf F}$ in $\FRC$) and parameters $p$ and $q$. 

Next, define $\WRC_{p,q}^{G, i}$ as the probability measure on $\{0,1\}^E$
corresponding to first setting $X(E \setminus E_i) \equiv 1$, and then
picking $X(E)$ in such a way that
\[
\WRC_{p,q}^{G, i}\bigl(X(E_i) = \xi\bigr) = 
\frac{q^{\kkk(\xi)}}{Z} \prod_{e\in E_i} p^{\xi(e)}(1-p)^{1-\xi(e)}
\]
where $\kkk(\xi)$ is the number of connected components of $\xi$ {\bf that
do not intersect $\partial V_i$}. 
Similarly as in (\ref{eq:monotone_limit}), we get
\[
\WRC_{p,q}^{G, 1} \geqd \WRC_{p,q}^{G, 2} \geqd \cdots
\]
(note the reverse inequalities), and thus also a limiting measure
$\WRC_{p,q}^G$ which we call the random-cluster measure on $G$ with
{\bf wired boundary condition} and parameters $p$ and $q$.

We now briefly discuss how the above relates to
the DLR approach to the random-cluster model on 
infinite graphs. It is
natural to expect that the limiting measures
$\FRC_{p,q}^G$ and $\WRC_{p,q}^G$ should satisfy some analogue of
(\ref{eq:single_edge_cond_prob}). Indeed, $\FRC_{p,q}^G$ admits
conditional probabilities satisfying
\begin{equation} \label{eq:DLR_free}
\FRC_{p,q}^G\bigl(X(e)=1 \bigm| X(E\setminus\{e\})= \xi\bigr) =
\left\{
\begin{array}{ll}
p & \mbox{if } x\leftrightarrow y, \\
\frac{p}{p+(1-p)q} & \mbox{otherwise}
\end{array} \right.
\end{equation}
for any $e\in E$ and any $\xi\in \{0,1\}^{E\setminus\{e\}}$, where the
event $x\leftrightarrow y$ is defined as in
(\ref{eq:single_edge_cond_prob}). 
Although $\WRC_{p,q}^G$ does
{\em not}, in general, satisfy the same local specification,
it satisfies
\begin{equation} \label{eq:DLR_wired}
\WRC_{p,q}^G\bigl(X(e)=1 \bigm| X(E\setminus\{e\})= \xi\bigr) =
\left\{
\begin{array}{ll}
p & \mbox{if } x\Ccon y, \\
\frac{p}{p+(1-p)q} & \mbox{otherwise,}
\end{array} \right.
\end{equation}
where $x\Ccon y$ denotes the event that either $\xi$ contains a path from 
$x$ to $y$, or it contains an infinite self-avoiding path starting at $x$
and an infinite self-avoiding path starting at $y$. In other words, 
$x\Ccon y$ is the same event as $x\leftrightarrow y$, except that 
in $x\Ccon y$ the path from $x$ to $y$ is allowed to go ``via infinity''. 
We think of this as a kind of compactification of the graph.
These facts are stated in \cite[Theorem 6.17]{GHM}.
(That $\FRC^G_{p,q}$ satisfies (\ref{eq:DLR_free}) is due to
\cite{BC:covariance}. The fact that $\WRC^G_{p,q}$ satisfies
(\ref{eq:DLR_wired}) can be proved analogously. Other proofs of \rref
e.DLR_wired/ can be found
in \cite{J} and in \cite{GHM}.) We 
call a probability measure on
$\{0,1\}^E$ a {\bf DLR random-cluster measure} (resp., a 
{\bf DLR wired-random-cluster measure}) 
with the given parameters $p$ and $q$
if it satisfies the local specifications in (\ref{eq:DLR_free})
(resp., in (\ref{eq:DLR_wired})).
(These local specifications are usually given on any finite
edge-set, rather than on a single edge. However, single-edge
specifications are enough; see, e.g., \cite[Theorem~6.18]{GHM}.) 
It turns out that $\FRC_{p,q}^G$ and
$\WRC_{p,q}^G$ play the following special role in the class of
DLR random-cluster and wired-random-cluster measures: If $\mu$ is any
DLR random-cluster measure or DLR wired-random-cluster measure for $G$
with parameters $p$ and $q$, then
\begin{equation} \label{eq:stoch_dom_1} 
\FRC_{p,q}^G \leqd \mu \leqd \WRC_{p,q}^G \, .
\end{equation}
We mention that (provided $G$ is connected)
the specifications (\ref{eq:DLR_free}) and 
(\ref{eq:DLR_wired}) differ with positive probability if and only if the
event
of having more than one infinite connected component has positive
probability.
By an application of the uniqueness theorem of Burton and Keane
\cite{BK}, we get in the case where $G$ is the usual $\Z^d$ lattice (and
more generally when $G$ is a transitive amenable graph --- see Section
\ref{sect:background_amenability_etc} for a definition) that the number
of infinite clusters is at most one, $\FRC$-a.s.\ as well as $\WRC$-a.s.
It follows that in this case, both $\FRC$ and $\WRC$ are simultaneously 
DLR random-cluster measures and DLR wired-random-cluster measures.

We finally state a few more well-known stochastic inequalities for
random-cluster measures that we shall have occasion to use:
If $p_1 \leq p_2$ and $p_1/[(1-p_1)q_1] \leq p_2/[(1-p_2)q_2]$, then
\begin{equation} \label{eq:witness_genfree}
\FRC^G_{p_1,q_1} \leqd \FRC^G_{p_2,q_2}
\end{equation}
and
\begin{equation} \label{eq:witness_genwired}
\WRC^G_{p_1,q_1} \leqd \WRC^G_{p_2,q_2}\,.
\end{equation}
In particular, if $p_1 \leq p_2$, then
\begin{equation} \label{eq:stoch_dom_3}
\FRC^G_{p_1,q} \leqd \FRC^G_{p_2,q}
\end{equation}
and
\begin{equation} \label{eq:stoch_dom_5}
\WRC^G_{p_1,q} \leqd \WRC^G_{p_2,q}\,.
\end{equation}

\subsection{The Potts model on infinite graphs}
\label{sect:background_infinite_Potts}

Let $G=(V,E)$ be infinite and locally finite, and
let $\{G_i:=(V_i, E_i)\}_{i=1}^\infty$ be as in the previous subsection. 
For $q\in \{2,3,\ldots\}$ and $\beta \geq 0$, 
define probability measures $\left\{\FPt_{q,
\beta}^{G,i}\right\}_{i=1}^\infty$
on $\{1,\ldots, q\}^V$ in such a way that the projection of
$\FPt_{q, \beta}^{G,i}$ on $\{1,\ldots, q\}^{V_i}$ equals
$\Pt_{q,\beta}^{G_i}$, and the spins on
$V\setminus V_i$ are i.i.d.\ uniformly distributed on $\{1, \ldots, q\}$
and independent of the spins on $V_i$. 
It turns out that $\FPt_{q, \beta}^{G,i}$
has a limiting distribution $\FPt_{q, \beta}^G$ as $i \rightarrow \infty$. 

Also, for a fixed spin $r\in \{1, \ldots, q\}$, define 
$\WPt_{q, \beta, r}^{G,i}$ to be the distribution corresponding to
picking $X\in \{1, \ldots, q\}^{V}$ by letting $X(V\setminus V_i) \equiv r$,
and letting $X(V_i)$ be distributed according to 
$\Pt_{q,\beta}^{G_i}$ {\em conditioned on the event that
$X(\partial V_i)\equiv r$}. Again, $\WPt_{q, \beta, r}^{G,i}$ has a limiting
distribution as $i\rightarrow \infty$, which we denote
$\WPt_{q, \beta,r}^G$. 

The existence of the limiting distributions
$\FPt_{q, \beta}^G$ and $\WPt_{q, \beta,r}^G$ are nontrivial results,
and in fact the shortest route to proving them goes via the stochastic 
monotonicity arguments for the random-cluster model outlined in Section
\ref{sect:background_infinite} and then using Propositions 
\ref{prop:from_FRC_to_FPotts} and \ref{prop:from_WRC_to_WPotts} below. 

A probability measure $\mu$
on $\{1, \ldots, q\}^V$ is said to be a Gibbs measure
(in the DLR sense) for the $q$-state Potts model on $G$ at inverse 
temperature $\beta$, if it admits conditional distributions such that
for all $v\in V$, all $r\in \{1, \ldots, q\}$, and all 
$\omega\in \{1, \ldots, q\}^{V\setminus \{v\}}$, we have 
\rlabel DLR_Potts
{\mu\bigl(X(v)=r \bigm| X(V\setminus \{v\}) = \omega\bigr) = \frac{1}{Z}
\exp\Bigl( -2\beta \sum_{\ed{ v,y } \in E}
\I{\omega(v) \neq \omega(y)} \Bigr)\,,
}
where the normalizing constant $Z$ may depend on $v$ and $\omega$ but not 
on $r$. It turns out that $\FPt_{q, \beta}^G$ and $\WPt_{q, \beta, r}^G$ 
are both Gibbs measures in this sense. Another Gibbs measure of particular 
interest is
\[
\WPt^G_{q, \beta} := \frac{1}{q} \sum_{r=1}^q \WPt_{q, \beta, r}^G \, . 
\]

The following extensions of Proposition \ref{prop:from_RC_to_Potts}
provide the relations between $\FRC$ and $\WRC$ on one hand, and
$\FPt$ and $\WPt$ on the other. 
\begin{prop} \label{prop:from_FRC_to_FPotts}
Fix an infinite locally finite graph $G$, an integer $q\geq 2$, and 
$p\in [0,1]$. Pick
a random edge configuration $X\in \{0,1\}^E$ according to $\FRC^G_{p,q}$. 
Then, for each connected component $\C$ of $X$, 
pick a spin uniformly from $\{1,\ldots, q\}$, and assign this spin to 
all vertices of ${\cal C}$. Do this independently for different
connected components. The $\{1, \ldots, q\}^V$-valued random spin
configuration arising from this procedure is then distributed according
to the Gibbs measure $\FPt^G_{q,\beta}$ for the $q$-state Potts model on $G$
at inverse temperature $\beta:= -\frac{1}{2}\log(1-p)$.  
\end{prop}
\begin{prop} \label{prop:from_WRC_to_WPotts}
Let $G$, $p$, $q$, and $\beta$
be as in Proposition \ref{prop:from_FRC_to_FPotts}.
Pick
a random edge configuration $X\in \{0,1\}^E$ according to the random-cluster
measure $\WRC^G_{p,q}$. Then, for each {\bf finite} 
connected component $\C$ of $X$, 
pick a spin uniformly from $\{1,\ldots, q\}$, and assign this spin to 
all vertices of ${\cal C}$. Do this independently for different
connected components. Finally assign value $r$ to all vertices of
infinite connected components. The $\{1, \ldots, q\}^V$-valued random spin
configuration arising from this procedure is then distributed according
to the Gibbs measure $\WPt^G_{q,\beta, r}$ for the $q$-state Potts model on
$G$
at inverse temperature $\beta$.
\end{prop}

\subsection{Some classes of infinite graphs}
\label{sect:background_amenability_etc}

The class of all infinite locally finite graphs is often too large to obtain
the most interesting results for the random-cluster model (and other 
stochastic models on graphs; see, e.g., \cite{L} for a survey), and indeed
most of our results will concern more restrictive classes of graphs.
Here we recall some such classes. 

Let, as usual, $G=(V,E)$ be an infinite locally finite graph. 
The number of edges incident to a vertex $x$ is called the {\bf degree} of
$x$. The graph
$G$ is said to be {\bf regular} if every vertex has the same degree.

A {\bf graph automorphism} of $G$ is a bijective mapping 
$\gamma: \, V \rightarrow V$ with the property that for all $x,y \in V$,
we have $\ed{ \gamma x, \gamma y } \in E$ if and only if
$\ed{ x,y } \in E$. Write $\Aut(G)$ for the group of all
graph automorphisms of $G$. To each $\gamma\in \Aut(G)$, there is a
corresponding mapping $\tilde{\gamma}: E\rightarrow E$ defined
by $\tilde{\gamma}\ed{ x,y }:= \ed{ \gamma x, \gamma y }$.
The graph $G$ is said to be {\bf transitive}
if and only if for any $x, y \in V$ there exists $\gamma\in\Aut(G)$
such that $\gamma x = y$. More generally, $G$ is said to be 
{\bf quasi-transitive} if $V$ can be partitioned into finitely many sets
$V_1, \ldots, V_k$ such that for any $i \in \{1, \ldots, k\}$ and any
$x,y \in V_i$, there exists a $\gamma\in\Aut(G)$
such that $\gamma x = y$.

A probability
measure $\mu$ on $\{0,1\}^E$ is said to be {\bf automorphism invariant} if
for
any $n$, any $e_1, \ldots, e_n \in E$, any $i_1, \ldots, i_n \in \{0,1\}$,
and any graph automorphism $\gamma$, we have
\[
\mu\bigl(X(e_1)=i_1, \ldots, X(e_n)=i_n\bigr) = 
\mu\bigl(X(\tilde{\gamma}(e_1))=i_1, \ldots,
X(\tilde{\gamma}(e_n))=i_n\bigr) \, .
\]
It follows from the construction of the free and wired random-cluster
measures
$\FRC$ and $\WRC$ (in particular from the
independence of the choice of $\{G_i=(V_i, E_i)\}_{i=1}^\infty$) that both
measures are automorphism invariant. It turns out that automorphism
invariance has far-reaching consequences for percolation processes on
various classes of transitive graphs; see, e.g., \cite{BK}, \cite{H3},
\cite{BLPS}, and \cite{LS}. 

Two important properties, that may or may not hold for a given quasi-transitive 
graph, are amenability and unimodularity, which we review next. We say that
a graph $G$ is {\bf amenable} if
\[
\inf \frac{|\partial W|}{|W|}=0 
\,,
\]
where the infimum ranges over all finite $W\subset V$ and $| \cdot |$
denotes cardinality. There are
various alternative definitions of amenability of a graph, which coincide
for transitive graphs (and more generally for graphs of bounded degree), but
not in general. 

For any graph $G$ and $x\in V$, define the {\bf stabilizer} $S(x)$ as the
set
of graph automorphisms that fix $x$, i.e.,
\[
S(x) := \{\gamma \in \Aut(G)  \st \gamma x = x \} \, .
\]
For $x,y \in V$, define
\[
S(x)y := \{ z \in V  \st \exists \gamma \in S(x) \mbox{ such that }
\gamma y = z \} \, .
\]
When $\Aut(G)$ is given the weak topology generated by its action on $V$, all
stabilizers are compact because $G$ is locally finite and connected.
We say that $G$ is {\bf unimodular} if for all $x,y$ in the same orbit of
$\Aut(G)$ (in the transitive case, this just means for all $x, y\in V$), we
have the symmetry
\[
|S(x)y| = |S(y)x| \, .
\]

An important class of transitive graphs is the class of
Cayley graphs. If $\Gamma$ is a finitely generated group with 
generating set
$\{g_1, \ldots, g_n\}$, then the {\bf Cayley graph} associated
with $\Gamma$ and that particular set of generators is the (unoriented)
graph
$G=(V,E)$ with vertex set $V:=\Gamma$, and edge set
\[
E:= \{ \ed{ x,y } \st x,y \in \Gamma, \exists i \in 
\{1, \ldots, n\} \mbox{ such that } x g_i = y \} \, .
\]
Most examples of graphs that have been studied in percolation theory
are Cayley graphs. These include $\Z^d$ (which, with a slight abuse of 
notation, is short for the graph with vertex set $\Z^d$ and edges connecting
pairs of vertices at Euclidean distance $1$ from each other), and the
regular 
tree $\T_n$ in which every vertex has exactly $n+1$ neighbours. The graph
$\Z^d$
is amenable, while $\T_n$ is nonamenable for $n\geq 2$. Also studied are
certain tilings of the hyperbolic plane 
(see Section \ref{sect:example}), and further examples can be
obtained, e.g., by taking Cartesian products of other Cayley graphs 
(such as $\T_n \times \Z$, the much-studied example of 
Grimmett and Newman \cite{GN}). 

All Cayley graphs are unimodular. An example, due to Trofimov \cite{Tr}, of
a transitive graph that is nonunimodular (and hence not a Cayley graph)
may be obtained by taking the binary tree $\T_2$, fixing a so-called 
topological end $\zeta$ (loosely speaking, a direction to infinity in the
tree),
and adding an edge between each vertex and its $\zeta$-grandparent.  

\section{The four critical values} \label{sect:critical}

Let, as usual, $G=(V,E)$ be infinite and locally finite. A probability
measure
$\mu$ on $\{0,1\}^E$ is said to be
{\bf insertion tolerant} if for any $e\in E$ and
almost every $\xi \in \{0,1\}^{E \setminus \{e\}}$, the conditional
$\mu$-probability that $e$ is open given the configuration $\xi$ on
$\{0,1\}^{E \setminus \{e\}}$, is strictly positive. Newman
and Schulman \cite{NS} showed that for any auto\-morph\-ism-invariant
insertion-tolerant
percolation process on $\Z^d$, the number of infinite clusters is
a.s.\ either $0$, $1$ or $\infty$. It has been observed by several authors
(see, e.g., \cite{BS}) that this result (as well as its proof in \cite{NS})
extends to the class of quasi-transitive connected graphs. 

Suppose that $G$ is quasi-transitive and connected. The Newman-Schulman result
then applies to $\FRC^G_{p,q}$ and $\WRC^G_{p,q}$ because
(\ref{eq:DLR_free}) and (\ref{eq:DLR_wired}) imply that the
two measures are insertion tolerant whenever $p\in (0,1)$. Furthermore,
$\FRC^G_{p,q}$ and $\WRC^G_{p,q}$ are ergodic; this 
was proved in \cite{BC:covariance} for $\FRC$ and in 
\cite{BBCK} for both measures. A simple proof of the stronger property
that the tail $\sigma$-field is trivial appears in \cite{L}. 
Ergodicity shows that for each fixed $p$ and $q$, the number
of infinite clusters is an
a.s.\ constant (which, however, need not be the same for $\FRC^G_{p,q}$ and 
for $\WRC^G_{p,q}$). Hence, given $G$ and $q$, the set $[0,1]$ of 
possible values
for $p$ can be partitioned into three sets according to whether the
$\FRC^G_{p,q}$-a.s.\ number of infinite clusters is $0$, $1$ or $\infty$,
and similarly it can be partitioned via $\WRC^G_{p,q}$.
From (\ref{eq:stoch_dom_3}) and (\ref{eq:stoch_dom_5}) we can immediately
deduce that the set of $p$ for which the number of infinite clusters is
$0$ is an interval containing $0$. In other words, there exist critical
values $\pc^\f:=\pc^\f(G,q)$ and $\pc^\w:=\pc^\w(G,q)$ such that
\begin{equation} \label{eq:p_c_free}
\FRC^G_{p,q} (\exists \mbox{ at least one infinite cluster})= \left\{
\begin{array}{ll}
0 & \mbox{for } p<\pc^\f, \\
1 & \mbox{for } p>\pc^\f
\end{array} \right. 
\end{equation}
and
\begin{equation} \label{eq:p_c_wired}
\WRC^G_{p,q} (\exists \mbox{ at least one infinite cluster})= \left\{
\begin{array}{ll}
0 & \mbox{for } p<\pc^\w, \\
1 & \mbox{for } p>\pc^\w \, .
\end{array} \right. 
\end{equation}
The question of how the interval above $\pc^\f$ ($\pc^\w$) is split
up according to whether the number of infinite clusters is $1$ or
$\infty$ is more intricate. Does it split nicely into two intervals, 
or are the sets more complicated? A proof is given in \cite[Proposition
5.2]{L} that in
the quasi-transitive unimodular case, the uniqueness sets are simply intervals.
Presumably, this holds whenever $G$ is quasi-transitive.
A similar proof shows that 
\rlabel{free_less_wired}
{
\FRC^G_{p,q}(\exists \mbox{ a unique infinite cluster}) \le
\WRC^G_{p,q}(\exists \mbox{ a unique infinite cluster})
\,.
}
Thus, if $G$ is quasi-transitive and unimodular, there
exist critical values $\pu^\f:=\pu^\f(G,q)$ and $\pu^\w:=\pu^\w(G,q)$ such
that
\begin{equation} \label{eq:p_u_free}
\FRC^G_{p,q} (\exists \mbox{ \rm a unique infinite cluster})= \left\{
\begin{array}{ll}
0 & \mbox{for } p<\pu^\f, \\
1 & \mbox{for } p>\pu^\f
\end{array} \right. 
\end{equation}
and
\begin{equation} \label{eq:p_u_wired}
\WRC^G_{p,q} (\exists \mbox{ \rm a unique infinite cluster})= \left\{
\begin{array}{ll}
0 & \mbox{for } p<\pu^\w, \\
1 & \mbox{for } p>\pu^\w \, .
\end{array} \right. 
\end{equation}

Summarizing (\ref{eq:p_c_free}), (\ref{eq:p_c_wired}), \rref e.p_u_free/,
and \rref e.p_u_wired/, we have (for $G$
connected, quasi-transitive and unimodular) four critical values
$\pc^\w, \pc^\f, \pu^\w, \pu^\f \in [0,1]$ such that
the $\FRC^G_{p,q}$-a.s.\ number of infinite clusters equals
\[
\left\{
\begin{array}{ll}
0 & \mbox{for } p<\pc^\f, \\
\infty & \mbox{for } p \in (\pc^\f, \pu^\f), \\
1 & \mbox{for } p> \pu^\f
\end{array} \right. 
\]
and the $\WRC^G_{p,q}$-a.s.\ number of infinite clusters equals
\[
\left\{
\begin{array}{ll}
0 & \mbox{for } p<\pc^\w, \\
\infty & \mbox{for } p \in (\pc^\w, \pu^\w), \\
1 & \mbox{for } p> \pu^\w \, .
\end{array} \right. 
\]

We now consider whether there is an
infinite cluster at $\pc(q)$. It is known that there is none for $q=1$ on
nonamenable quasi-transitive unimodular graphs \cite{BLPS, BLPScrit}. On the
other hand, it is known that there can be infinitely many
infinite clusters for $q > 2$
on the Cayley graphs $\T_n$ for $n \ge 2$
with respect to the wired random-cluster measure; see \cite{CCST,H2}.
While we do not have a criterion that settles the question completely, 
we have the following partial results:

\procl t.death
Let $G$ be a quasi-transitive nonamenable unimodular graph and $q \ge 1$. 
Then there is no infinite cluster $\FRC^G_{\pcf(q),q}$-a.s.
Also, the following are equivalent:
\begin{itemize}
\item[(i)] There is no infinite cluster $\WRC^G_{\pcw(q),q}$-a.s.
\item[(ii)] For every edge $e \in G$, the function $p \mapsto
\WRC^G_{p,q}\big(X(e)=1\big)$ is continuous from the left at $\pcw(q)$.
\item[(iii)] $\FRC^G_{\pcw(q),q} = \WRC^G_{\pcw(q),q}$.
\end{itemize}
\endprocl

\proof
Fix an edge $e = \ed{x,y}$.
Let $Q^\f(p) := \FRC^G_{p,q}\big(X(e)=1\big)$ and
$Q^\w(p) := \WRC^G_{p,q}\big(X(e)=1\big)$.
Now 
\[
Q^\f(p) = \lim_{i \to\infty} \FRC^{G, i}_{p,q}\big(X(e)=1\big)
\] 
by definition.
Since the latter probabilities (before the limit is taken)
are rational functions in
$p$ and $q$, they are continuous at all $p$. Since they are increasing in
$p$ and increase to their limit, we obtain left-continuity of
$Q^\f(p)$
at all $p$.

We need the existence of an invariant coupling
$(X,X')$ that witnesses the stochastic domination
$\FRC_{\pc^\f (q),q} \geqd \FRC_{p,q}$, and similarly for $\WRC$.
Such couplings are constructed in \cite{HJL}.
For $p < \pcf(q)$, the probability that 
\rlabel disagree
{\Xf{\pc^\f(q)}(e) = 1 \mbox{ and } \Xf{p}(e) = 0
}
in any such coupling equals $Q^\f\big(\pcf(q)\big)-Q^\f(p)$. Thus,
left-continuity of $Q^\f(p)$ at $p = \pc^\f(q)$ implies that
the probability of \rref e.disagree/
tends to 0 as $p \uparrow \pc^\f(q)$. Therefore, minor
modifications of the proofs in \cite{BLPScrit} show that there is no
infinite cluster $\FRC^G_{\pcf(q),q}$-a.s.
This establishes the first claim.

Now we consider the equivalences.
A similar argument to the above shows that (ii) implies (i). 
The other implications hold for all graphs.
That (i) implies (iii) is due to \cite{ACCN}; the reasoning is
sketched in the proof of \rref p.f=w/ below.
If (iii) holds, then $\FRC_{p,q}=\WRC_{p,q}$ for all $p \le \pcw$
by \cite{ACCN} again. Therefore $Q^\w(p) = Q^\f(p)$ for $p \le \pcw$ and
the continuity of $Q^\f(p)$ implies (ii) as above.
%
\Qed

Given a spin configuration $\omega \in\{1, \ldots, q\}^V$ and an edge 
configuration $\xi \in \{0,1\}^E$, we may partition $V$ into {\bf
connected
single-spin components}, meaning that $x$ and $y$ are in the same 
connected single-spin component if and only if there is a path from
$x$ to $y$ in which all vertices have the same spin and all edges are
open. 
The following facts relate the four critical values to 
corresponding phenomena for the Potts model. 

\procl p.facts
Let $G$ be any graph,
$\beta > 0$, $q \ge 1$, and $p := 1-e^{-2\beta}$.
\begin{itemize}
\item[(i)] There is no infinite cluster 
$\WRC_{p,q}$-a.s.\ iff there is a unique Gibbs measure for the Potts model with
the corresponding parameters $(q, \beta)$.
\item[(ii)] Let $\omega \in \{1, \ldots, q\}^V$ be chosen according to $\FPt_{q,
\beta}$ and
independently $\xi \in \{0,1\}^E$ be chosen according to Bernoulli$(p)$
percolation. There is a unique infinite cluster $\FRC_{p,q}$-a.s.\
iff $(\omega,\xi)$ a.s.\ produces a unique infinite connected
single-spin component.
\item[(iii)] Let $\omega \in \{1, \ldots, q\}^V$ be chosen according to
$\WPt_{q, \beta}$ and
independently $\xi \in \{0,1\}^E$ be chosen according to Bernoulli$(p)$
percolation. There is a unique infinite cluster $\WRC_{p,q}$-a.s.\
iff $(\omega,\xi)$ a.s.\ produces a unique infinite connected
single-spin component.
\item[(iv)] If there is a unique infinite cluster $\FRC_{p,q}$-a.s., then
$\FPt_{q, \beta}$ is not extremal among all Gibbs measures.
\item[(v)] Suppose that $G$ is quasi-transitive.
If there is a unique infinite cluster $\FRC_{p,q}$-a.s., then
$\FPt_{q, \beta}$ is not extremal among invariant Gibbs measures.
If $G$ is also unimodular, then the converse holds; in fact, if there is not
a unique infinite cluster $\FRC_{p,q}$-a.s., then
$\FPt_{q, \beta}$ is ergodic, i.e., extremal among all invariant
probability measures.
\end{itemize}
\endprocl

\proof
Part (i) is essentially due to Aizenman et al.\ \cite{ACCN} (or see 
\cite{GHM}). 
Parts (ii) and (iii) follow immediately from the coupling of
random-cluster and Potts models underlying Propositions
\ref{prop:from_FRC_to_FPotts} and \ref{prop:from_WRC_to_WPotts}.

Part (iv) follows from (ii): If there is a unique infinite component
$\FRC_{p,q}$-a.s., then let $(\omega,\xi)$ have the distribution
$\FPt_{q,\beta}\times\RC_{p,1}$, as in (ii). By (ii), we may define
$r(\omega,\xi)$ to be the spin of the
unique infinite single-spin component determined by $(\omega,\xi)$.
Let $\CC(\omega)$ be the collection of
maximal connected subgraphs of $G$ whose vertices have a common spin; this
does not depend on $\xi$.
Then $r(\omega,\xi)$
is the spin of the unique graph $K$ in $\CC(\omega)$ such that $K \cap \xi$
contains an infinite component with positive probability (in $\xi)$.
In particular, $r(\omega,\xi)$ depends only on $\omega$ a.s.; and it is
a tail random variable that is not trivial. 
Thus, we have shown that the tail $\sigma$-field of $\FPt_{q,\beta}$ is not
trivial. This is equivalent to nonextremality among all Gibbs measures by
\cite[Theorem 7.7]{Geor}.

The first part of (v) is due to \cite[Theorem 4.2]{Schon:mult}.
The converse of (iv) is not true in general, as is well known on trees
(see, e.g., \cite{BRZ} or \cite{EKPS}). 
However, if $G$ is quasi-transitive and unimodular and
if there is not a unique infinite cluster $\FRC_{p,q}$-a.s., then by
\cite[Theorem 4.1 and Lemma 6.4]{LS} extended from the transitive case to
the quasi-transitive case, $\FPt_{q, \beta}$ is ergodic,
which is the same
as extremal among all invariant measures.
\Qed

How do the four critical values relate to each other? From the definitions
it is immediate that 
\begin{equation} \label{eq:inequality_f}
\pc^\f \leq \pu^\f
\end{equation}
and 
\begin{equation} \label{eq:inequality_w}
\pc^\w \leq \pu^\w \, .
\end{equation}
By (\ref{eq:stoch_dom_1}) and \rref
e.free_less_wired/, we also have
\begin{equation} \label{eq:inequality_c}
\pc^\w \leq \pc^\f
\end{equation}
and 
\begin{equation} \label{eq:inequality_u}
\pu^\w \leq \pu^\f \, .
\end{equation} 
All of (\ref{eq:inequality_f})--(\ref{eq:inequality_u})
can reduce to equalities; this happens, e.g., 
whenever $G$ is amenable. To see this for (\ref{eq:inequality_f}) and
(\ref{eq:inequality_w}), just note the well-known fact 
that the Burton-Keane \cite{BK,GKN} encounter-point
argument (for showing uniqueness of the infinite cluster under the
insertion tolerance condition) goes through in the amenable setting. For
(\ref{eq:inequality_c}) and (\ref{eq:inequality_u}), 
see Grimmett \cite{Gr} and Jonasson \cite{J}, where it is shown that for
all $q \ge 1$, there are at most countably many $p$ such that
$\FRC^G_{p,q} \ne \WRC^G_{p,q}$. 

The inequalities can also be strict. To get examples with strict
inequalities
in (\ref{eq:inequality_f}) and (\ref{eq:inequality_w}), one can take $q:=1$
and $G$ to be any of the nonamenable transitive unimodular graphs that are
known to have a ``middle phase'' for i.i.d.\ percolation (i.e., a
positive-length interval of
values of $p$ that give rise to infinitely many infinite clusters); see,
e.g.,
\cite{L}. Using the ideas in the proof of \cite[Proposition 5.2]{L}
and the inequalities \rref
e.witness_genfree/ and \rref e.witness_genwired/, it is not hard to show
that
one can take $q$ to be slightly larger than $1$ in all such examples. 
For an example where the inequality in (\ref{eq:inequality_c}) 
is strict, we can simply take $G$ to be the regular
tree $\T_n$ with $n\geq 2$ and $q>2$ (see, e.g., H\"aggstr\"om \cite{H2}),
or
take any nonamenable regular graph
with $q$ sufficiently large (this follows from \cite[Theorem 1.2(a)]{J}
in combination with (\ref{eq:stoch_dom_1})). 
Finally, for an example where (\ref{eq:inequality_u}) is strict, we
refer to Section \ref{sect:example}. 

The inequalities (\ref{eq:inequality_f})--(\ref{eq:inequality_u}) say
nothing about the relation between  $\pu^\w$ and $\pc^\f$. Here it is
possible to get a strict inequality in either direction. Examples with
$\pc^\f< \pu^\w$ can be obtained in the same way as for 
(\ref{eq:inequality_f}) and (\ref{eq:inequality_w}); see \rref p.small_q/
below for explicit bounds on certain graphs. To get an example
with the reverse inequality $\pu^\w < \pc^\f$ is more intricate and is the
topic of the next section.
Note that any such example also gives strict inequality in \rref
e.inequality_c/ and in \rref e.inequality_u/.

\procl q.boldconj
We say that $G$ has {\bf one end} if the complement of every finite
subset has exactly one infinite component.
If $G$ is any nonamenable quasi-transitive graph with one end and $q \ge
1$, are the inequalities {\rm \rref e.inequality_f/} and 
{\rm \rref e.inequality_w/}
necessarily strict?
\endprocl

Of course, when $q=1$, a famous conjecture of \cite{BS} asserts a
positive answer.

If $G$ is a graph drawn in the
plane in such a way that edges do not cross and such that each bounded set
in the plane contains only finitely many vertices of $G$, then $G$ is said
to be {\bf properly embedded}. We shall always assume without mention that
planar graphs are properly embedded. (The graphs we shall consider in the
next section can be
embedded in the hyperbolic plane more geometrically than in the euclidean
plane, but topologically and combinatorially, this is not different from
euclidean embeddings.)
If $G$ is a planar (multi)graph, then the
{\bf planar dual} $\d G$ of $G$ (really, of this particular embedding of
$G$) is the (multi)graph formed as follows: The vertices of $\d G$ are the
faces formed by $G$. Two faces of $\d G$ are joined by an edge precisely
when they share an edge in $G$. Thus, $E(G)$ and $E(\d G)$ are in a
natural one-to-one correspondence. Furthermore, if one draws each
vertex of $\d G$ in the interior of the corresponding face of $G$
and each edge of $\d G$ so that it crosses the corresponding edge of
$G$, then the dual of $\d G$ is $G$. 
For planar graphs, we shall always assume that $G$ and its planar dual $\d
G$ are locally finite, whence each graph has one end.

Let $G$ be a planar graph.
If $\mu$ is a probability measure on $\{0,1\}^E$, we associate a {\bf dual}
measure $\d \mu$ on $\{0,1\}^{\d E}$ as follows. Given $e \in E$, let $\d e$
be the edge in $\d E$ that crosses $e$.
Given $\xi \in \{0, 1\}^E$, let $\tilde \xi \in \{0, 1\}^{\d E}$ be the
function $\d e \mapsto 1 - \xi(e)$.
For a Borel set $A \subset \{0,1\}^E$, write $\tilde A := \{\tilde \xi
\st \xi \in A\}$.
Then $\d \mu$ is defined by $\mu(A) = \d\mu(\tilde A)$.
Our next proposition is more or less well known (see, e.g., \cite{CCS,Wel}),
but perhaps has not been stated in this particular form before. 
For completeness, we provide the simple proof here. 

\procl p.FWdual
For any planar graph $G$, $\FRC^G_{p,q}$ is dual to $\WRC^{\d G}_{p',q}$ if 
\rlabel p'
{p' = {(1-p)q \over p+(1-p)q}\,.
}
\endprocl

\proof
Suppose first that $G$ is a finite graph.
For any $\xi \in \{0,1\}^E$, write $\tilde\xi :=
\d E \setminus \d \xi$ as above. Euler's formula applied to the graph
$\tilde G := (\d V, \tilde\xi)$ says that 
$$
|\d V|-|\tilde \xi|+\kk(\xi) = 1 + \kk(\tilde \xi)
$$
since the number of faces of $\tilde G$ is equal to
$\kk(\xi)$. Thus,
\begin{align*}
\RC^G_{p,q}(\xi)
&=
Z^{-1} p^{|\xi|} (1-p)^{|E\setminus\xi|} q^{\kk(\xi)}
=
Z^{-1} p^{|\d E \setminus \tilde\xi|} (1-p)^{|\tilde\xi|} q^{1 + \kk(\tilde
\xi)
- |\d V| + |\tilde\xi|}
\cr&=
Z^{-1} q^{1 - |\d V|} p^{|\d E \setminus \tilde\xi|}
[(1-p)q]^{|\tilde\xi|} q^{\kk(\tilde \xi)}
=
\tilde Z^{-1} (1 - p')^{|\d E \setminus \tilde\xi|}
{p'}^{|\tilde\xi|} q^{\kk(\tilde \xi)}
=
\RC^{\d G}_{p',q}(\tilde \xi)
\,,
\end{align*}
where $Z$ and $\tilde Z$ are normalizing constants that do not depend on
$\xi$.
Thus, $\RC^G_{p,q}$ is dual to $\RC^{\d G}_{p', q}$.

Now let $V_i$ be a sequence of increasing finite subsets of $V$ such that
the faces of $G(V_i)$ are faces of $G$, except for the outer face, of
course. Then we have seen that $\RC^{G_i}_{p, q}$ is dual to
$\RC^{\d G_i}_{p', q}$. Thus, the general result follows by taking
weak limits.
\Qed

In consequence, the methods of Benjamini and Schramm \cite{BS:hp} show the
following proposition. 
The paper \cite{BS:hp} dealt only with transitive graphs, but
the methods work just as well given the fact that quasi-transitive planar
graphs with one end are unimodular; this fact is proved in \rref
b.LP:book/.

\procl p.FWcoupled
Let $G$ be a planar nonamenable 
quasi-transitive graph and $p'$ be as in
{\rm \rref
e.p'/}. In the natural coupling of $\FRC^G_{p,q}$ and $\WRC^{\d
G}_{p',q}$ as dual measures, the number of infinite clusters with respect
to each is a.s.\ one of the following: $(0,1)$, $(1,0)$, or
$(\infty,\infty)$.
\endprocl

Write 
\rlabel{hdef}
{
h(x) := x/(1-x)
\,.
}

\procl c.pcpu
For any planar nonamenable quasi-transitive graph $G$,
$$
h\bigl(\pcw(G, q)\bigr) h\bigl(\puf(\d G, q)\bigr) = 
h\bigl(\pcf(G, q)\bigr) h\bigl(\puw(\d G, q)\bigr) = q
\,,
$$
$0 < \pcw(G, q) \le \pcf(G, q) < 1$, and $0 < \puw(G, q) \le \puf(G, q) <
1$.
\endprocl

\proof
\rref p.FWcoupled/ shows that there is
no infinite cluster a.s.\ with respect to the free measure iff there
is a unique infinite cluster a.s.\ with respect to the dual wired measure.
Therefore, $\puw(\d G, q) = \left(\pcf(G, q)\right)'$ in the notation of
\rref e.p'/. Some algebra shows that this is the same as $h\bigl(\pcf(G,
q)\bigr) h\bigl(\puw(\d G, q)\bigr) = q$. A similar proof shows the other
equation.

It is well known that $0 < \pc(G) := \pc(G, 1) < 1$ under the present
assumptions
\cite{L}. From \rref e.witness_genfree/ and \rref e.witness_genwired/, it
follows that the same holds for $\pcw(G, q)$ and $\pcf(G, q)$ when $q > 1$.
The same now follows for the uniqueness critical
points by the equations just established.
\Qed

\procl c.uniq_p_u
Let $G$ be a planar nonamenable quasi-transitive graph with one end.
Then there is
$\WRC^G_{\puw(G, q),q}$-a.s.\ a unique infinite cluster.
\endprocl

\proof
In light of \rref t.death/, there is no infinite cluster $\FRC^{\d
G}_{\pcf(\d G, q),q}$-a.s.
Hence by \rref p.FWcoupled/ and \rref c.pcpu/, there is a unique infinite
cluster $\WRC^G_{\puw(G, q),q}$-a.s.
\Qed

The methods of 
Grimmett \cite{Gr:compare} show that for any quasi-transitive $G$, the
critical points $\pcw(q)$ and $\pcf(q)$ are 
continuous and strictly increasing functions of $q$ as long as
$\pc(1) < 1$, and similarly $\puw(q)$ and $\puf(q)$ are 
continuous and strictly increasing functions of $q$ as long as $\pu(1) < 1$.
When $G$ is a planar regular graph with regular dual, \cite{BS:hp} shows
that $\pc(G) < \pu(G) := \pu(G, 1)$, while \rref t.separated/ below shows that
$\pcf(G, q) > \puw(G, q)$ for large $q$. Thus, there is at least one $q$
for which $\pcf(G, q) = \puw(G, q)$. If $G$ is isomorphic to its dual,
then these critical values are equal to $\sqrt q/(\sqrt q + 1)$ because
of \rref c.pcpu/. We do not know whether this holds for exactly one $q$; see
\rref q.pupcphase/.

Another natural set of $p$ to examine for fixed $q$ is the set where
$\FRC_{p, q} = \WRC_{p, q}$, which is where
$\FPt_{q, \beta} = \WPt_{q, \beta}$, as shown in
\cite{J}, Lemma 4.3. In general, this is not simply an interval, but
let us define
$$
\peq(G) := \sup \{ p < 1 \st \FRC^G_p \ne \WRC^G_p \}\,.
$$
We say that a graph has {\bf bounded fundamental cycle length} if
there is a set of (oriented simple) cycles of the graph with bounded
length that generates all cycles by addition and subtraction (in the sense
of homology).
For example, this is the case for Cayley graphs of finite
presentations of groups.
Recall that $G$ has one end if the complement of every finite
subset has exactly one infinite component.

\procl p.f=w
Let $G$ be any graph. 
\begin{itemize}
\item[(i)] If $p < \pcw(G)$, then $\FRC^G_p = \WRC^G_p$.
\item[(ii)] If $\pcw(G) < p < \puf(G)$, then $\FRC^G_p \ne \WRC^G_p$. 
\item[(iii)] If $G$ is a graph with one end and bounded fundamental cycle
length, then $\peq(G) < 1$. 
\item[(iv)] If $G$ is a planar nonamenable 
quasi-transitive graph with one
end, then $\peq(G) = \puf(G)$.
\end{itemize}
\endprocl

\proof Part (i) is due to \cite{ACCN}, but we recap the short
proof here. It suffices to show that $\WRC^G_p \leqd \FRC^G_p$ if $p <
\pcw(G)$, or, more generally, if there is no infinite cluster
$\WRC^G_p$-a.s.  Since there is no infinite cluster, given any ball $B$
about $o$ and any $\epsilon > 0$, there is a ball $B'$ so that with
probability at least $1 - \epsilon$, there is a unique maximal set $K
\subset B'$ such that all of $\bde K$ is closed and $B \subset K$. Given
that $K$ is such a set, the configuration restricted to $G(K)$ has the
distribution $\RC^{G(K)}_p$, which is dominated by the restriction of
$\FRC^G_p$ to $G(K)$.  In particular, this holds for the restriction of
the configuration to $B$.  Since $\epsilon$ and $B$ were arbitrary, the result
follows.

Part (ii) is
due to \cite{J}. Again, the proof is short: If $\FRC^G_p = \WRC^G_p$,
then the two measures give the same number of infinite clusters a.s.
Furthermore, since the measures are both DLR random-cluster measures
and DLR wired random-cluster measures, there is at most one infinite
cluster a.s. Hence $p \notin \bigl(\pcw(G), \puf(G)\bigr)$.

Now let $G$ be any graph with one end and bounded fundamental cycle
length. 
Let $t$ be an upper bound for the lengths of a set of generating cycles.
The fundamental theorem of \cite{BB} implies that the $t$-neighborhood of any
minimal cutset is connected.
Note that by \rref e.witness_genfree/, $\FRC^G_p$
dominates Bernoulli($p^*$) bond percolation $\P_{p^*}$ with $p^*$ close to
1 when $p$ is close to 1. 
Consider then a coupling $(\xi, \omega)$ with $\xi \supseteq
\omega$ such that $\xi$ has distribution
$\FRC^G_p$ and $\omega$ has distribution $\P_{p^*}$.
Let $\omega'$ consist of those edges $e$ such that all bonds
within distance $t$ of $e$ are open in $\omega$; define $\xi'$ similarly
with respect to $\xi$.
The percolation $\omega'$ dominates a Bernoulli percolation with
survival parameter close to 1 by \cite{LSS} or \cite[Remark 6.2]{LS}.
Thus, a.s.\ the {\it closed\/} bonds of $\omega'$ do not form any infinite
cluster if $p^*$ is sufficiently close to 1.
Therefore, when $p$ is close to 1, the same holds for $\xi'$.
Choose $p$ so close to 1 that $\xi'$ has no infinite closed clusters.
Fix $r \in \Z^+$ and $\epsilon > 0$.
Let $B_r$ be the ball of radius $r$ about $o$.
For any set of sites $S$, let $K(S)$ be 
the set of vertices that can be reached from some vertex
of $S$ without using an open bond from $\xi'$. 
There is a sphere $S$ about $o$ so that the probability of the event $E :=
\{K(S) \cap B_{r+t} = \emptyset\}$ is at least $1 - \epsilon$.
Since $S$ separates $B_r$ from $\infty$ as a set of vertices, so does $\bde
S$ as a set of edges.
Hence the same holds for the larger set $K(S)$ on the event $E$.
Let $L$ be the $t$-neighborhood of $\bde K(S)$.
Then on the event $E$, the set $L$ consists of open edges of $\xi$ and $L$
separates $B_r$ from $\infty$. 
Furthermore, some minimal cutset of $\bde K(S)$ separates $B_r$ from
$\infty$ and has a connected $t$-neighborhood, whence $L$ contains a
connected open cutset that separates $B_r$ from $\infty$.
Therefore, given $E$, the configuration $\xi$ restricted to $G(B_r)$ has the
distribution $\WRC^{G,B_r}_p$ [which means $\WRC^{G,i}_p$ if $G_i = G(B_r)$],
which dominates the restriction of $\WRC^G_p$ to $G(B_r)$.  Since
$\epsilon$ and $r$ were arbitrary, part (iii) follows.

Finally, if $G$ is planar, nonamenable, and quasi-transitive, consider $p >
\puf(G)$. Let $p'$ be as in \rref e.p'/.
We have $p' < \pcw(\d G)$ by \rref c.pcpu/, whence $\FRC^{\d
G}_{p'} = \WRC^{\d G}_{p'}$ by part (i). Because
of \rref p.FWdual/, it follows that $\WRC^G_p = \FRC^G_p$.
Hence $\peq(G) = \puf(G)$ by part (ii).
\Qed

\procl r.conj_Schon
The argument of the last paragraph shows that as long as there is a unique
infinite cluster $\FRC_{\puf}$-a.s., then $\WRC^G_p = \FRC^G_p$, which
establishes Conjecture 4.2 of \cite{Schon:mult} in the planar
nonamenable case.
\endprocl

\procl q.FequalsW
If $G$ is a quasi-transitive graph and $q > 1$, is the set of $p \in [0,1]$
where $\FRC^G_p \ne \WRC^G_p$ an interval?
\endprocl

\section{Isoperimetric constants and the critical values}
\label{sect:example}

In order to show that there is a Cayley graph $G$ with $\pcf(q) >
\puw(q)$ for some $q > 1$, we shall need an estimate of an
isoperimetric constant. In fact, we are able to calculate precisely
the necessary isoperimetric constants for planar regular graphs whose
dual is also regular (in this case, either the graph or its dual is a
Cayley graph; see \cite{CK}). Planar duality and Euler's formula will
be essential for this. 

We shall make use of the following isoperimetric constants.
For $K \subseteq V$, recall that
$E(K) := \{\ed{x,y} \in E \st x, y \in K\}$ and set
$E^*(K) := \{\ed{x, y} \in E \st x \in K \hbox{ or } y \in K\}$.
Define $\bde K := E^*(K) \setminus E(K)$
and $G(K) := \bigl(K, E(K)\bigr)$.
Write
\[
\isoe'(G) := \lim_{N \to\infty} \inf \left\{ \frac{|\bde K|}{|K|} \st
  K \subset V,\ G(K) \hbox{ connected}, \ N \le |K| < \infty \right\}\,,
\]
\[
\beta(G) := \lim_{N \to\infty} \inf \left\{ \frac{|K|}{|E(K)|} \st
  K \subset V,\ G(K) \hbox{ connected}, \ N \le |K| < \infty \right\}\,,
\]
\[
\delta(G) := \lim_{N \to\infty} \sup \left\{ \frac{|K|}{|E^*(K)|} \st
  K \subset V,\ G(K) \hbox{ connected}, \ N \le |K| < \infty \right\}\,.
\]
We write $d_G$ for the degree of vertices in $G$ when $G$ is regular.
We remark that $\beta(G) = 2/\alpha(G)$, with $\alpha(G)$ defined as in
\cite{BLPS}, except that $\alpha$ was defined with an infimum, rather
than a liminf. In any case, when $G$ is regular, 
\rlabel betaisoe
{\beta(G) = {2 \over d_G - \isoe'(G)}
}
and
\rlabel deltaisoe
{\delta(G) = {2 \over d_G + \isoe'(G)}\,.
}
It is shown in \cite{BLPS} that when $G$ is transitive, 
\[
\isoe'(G) = \inf \left\{ {|\bde K| \over |K|} \st K \subset V \hbox{
finite and nonempty} \right\}\,.
\]
(The right-hand side is denoted $\isoe(G)$ there.) Thus, when $G$ is
transitive, we have that
\rlabel deltasup
{\delta(G) = \sup \left\{ {|K| \over |E^*(K)|} \st K \subset V \hbox{
finite and nonempty} \right\}\,.
}
Recall from Section 
\ref{sect:background_amenability_etc} that $G$ is called 
quasi-transitive if the vertex set of $G$ decomposes
into a finite number of orbits under the action of $\Aut(G)$. Note
that $G$ is quasi-transitive iff $\d G$ is quasi-transitive. 

The estimate that we shall need is embodied in \rref c.beta_beta/, but
the precise combinatorial calculation is the following.

\procl t.isocalc
If $G$ is a planar regular graph with regular dual $\d G$, then
$$
\isoe'(G) = (d_G-2) \sqrt{1 - {4 \over (d_G - 2)(d_{\d G} - 2)}}\,.
$$
\endprocl

\procl r.reg_coreg 
In this case, $G$ and $\d G$ are transitive. 
This is folklore. Since we have been unable to find a suitable reference,
we include a proof here.
First, recall the existence of tessellations by {\it congruent} polygons
(in the euclidean or hyperbolic plane, as necessary).  It is
easy to see that the edge graphs of any two such tessellations of the same
type are isomorphic, by going out ring by ring around a starting polygon,
and thus that such edge graphs are transitive.
Now we assert that any (proper) tessellation of a plane with degree $d$ and
codegree $\d d$ has an edge graph that is isomorphic to the edge graph of
the corresponding tessellation above. In case $(d-2)(\d d-2)=4$, we replace
each face by a congruent copy of a flat polygon; in case $(d-2)(\d d-2)>4$,
replace it by a congruent copy of a regular hyperbolic polygon (with
curvature $-1$) of $\d d$ sides and interior angles $2\pi/d$;
while if $(d-2)(\d d-2)<4$,
replace it by a congruent copy of a regular spherical polygon (with
curvature $+1$) of $\d d$ sides and interior angles $2\pi/d$. 
Glue these together along the edges. We get a
Riemannian surface of curvature 0, $-1$, or $+1$, correspondingly, that is
homeomorphic to the plane since our assumption is that the plane is the
union of the faces, edges, and vertices of the tessellation, without
needing any limit points.
Riemann's theorem says that the surface is isometric to either
the euclidean plane or the hyperbolic plane (the spherical case is
impossible). That is, we now have a
tessellation by congruent polygons.
(One could also prove the existence statement in
a similar manner.)
\endprocl

\procl r.balls
Contrary to what one might first expect, we believe that combinatorial
balls never give the best isoperimetric constants when $G$ is a
nonamenable planar transitive graph with one end.
For example, for the isoperimetric constant $\isoe'(G)$, we believe that
$\liminf_n {|\bde B_{n}| / |B_n|} >
\isoe'(G)$, where $B_n$ is the ball of radius $n$ about $o$.
In many cases, this follows
from the formulas of \cite{FP}. For instance, suppose that
$$
\theta
:=
\lim_{n \to\infty} {|B_{n+1}| \over |B_n|}
$$
exists. (This is not the case for all Cayley graphs; see \cite{GH}.)
Then
$$
\theta - 1
\le \liminf_{n \to\infty} {|\bde B_{n}| / |B_n|}
\,.
$$
Thus, if 
$\theta - 1 > \isoe'(G)$, then we may conclude that
$\liminf_{n \to\infty} {|\bde B_{n}| / |B_n|} > \isoe'(G)$.
For example, if $d_{\d G} = 6$, in which case $G$ is always a Cayley graph
(see \cite{CK}), then
\rref t.isocalc/ shows that $\isoe'(G) = \sqrt{(d_G-2)(d_G-3)}$. On 
the other hand, \cite{FP} shows that
\begin{align*}
\sum_{n \ge 0} |B_{n}\setminus B_{n-1}| z^n
&=\frac{z^3+2z^2+2z+1}{z^3+(2-d_G)(z^2+z)+1}
=
\frac{z^2+z+1}{z^2+(1-d_G)z+1}
\\&=
\frac{z^2+z+1}{(1-\gamma z)(1-\gamma^{-1} z)}
=
(z^2+z+1) \sum_{n \ge 0} (\gamma z)^n \sum_{m \ge 0} (\gamma^{-1} z)^m
\,,
\end{align*}
where $\gamma$ is the smallest positive root of $z^2+(1-d_G)z+1$.
Therefore,
$$
|B_{n}\setminus B_{n-1}| =
\gamma^n(3-\gamma^{-2n-2}-\gamma^{-2n}-\gamma^{-2n+2})/(1-\gamma^{-2})
$$
for $n \ge 1$.
Thus, $\theta = \gamma$ exists and
$\theta - 1 = \sqrt{d_G - 3}\bigl(\sqrt{d_G - 3} +
\sqrt{d_G+1}\bigr)/2$, which is easily verified to be larger than
$\isoe'(G)$.
One may also verify from the lower bound of \cite[Theorem
5.1]{BauesPeyerimhoff} that if $\d d  - 2 \ge d \ge 5$, then $\inf_n |B_n
\setminus B_{n-1}|/|B_n| > \isoe(G)$. 
\endprocl

\rref t.isocalc/ follows from applying the following identity to $G$
and $\d G$, then solving the resulting two equations using \rref
e.betaisoe/ and \rref e.deltaisoe/:

\procl t.beta_delta
For any planar regular graph $G$ with regular dual, we have
$$
\beta(G) + \delta(\d G) = 1\,.
$$
\endprocl

\proof
Note first that the constant $\isoe'(G)$ is unchanged if, in its
definition, we require $K$ to be connected and simply connected when $K$ is
regarded as a union of closed faces of $\d G$ in the plane. This is because
filling in holes increases $|K|$ and decreases $|\bde K|$. Since $G$ is
regular, the same holds for $\beta(G)$ by \rref e.betaisoe/. Likewise, the
assumed regularity of $\d G$ and \rref e.deltaisoe/ imply a comparable
statement for $\delta(\d G)$. In fact, we shall need a refinement of this
idea for $\delta(\d G)$. Namely, given a finite connected set $K$ in $V(\d
G)$, regard each element of $K$ as a face of $G$ and
let $K' \subset V$ be the set of vertices bounding these faces.
Let $\widehat K$ be the set of all faces in $\d G$ that lie in the interior
of
the outermost cycle formed by $E(K')$. Then again $|\widehat K| \ge |K|$ and
$|\bde \widehat K| \le |\bde K|$, so that $\delta(\d G)$ can be approached
arbitrarily closely by such sets $\widehat K$. Note that $|E^*(\widehat K)|
=
|E((\widehat K)')|$.

Now let $\epsilon > 0$ and let $K$ be a finite connected set
in $V(\d G)$ such that $|E^*(K)| > 1/\epsilon$,
$|K|/|E^*(K)| \ge \delta(\d G) - \epsilon$, and 
$|E(K')| = |E^*(K)|$, where $K'$ is defined as above.
Since the number of faces of the graph $G(K')$ is
at least $|K|+1$, Euler's formula applied to the graph $G(K')$
gives
\rlabel upperbnd
{|K'|/|E(K')| + |K|/|E^*(K)| \le 1 + 1/|E^*(K)| < 1 + \epsilon\,.
}
Our choice of $K$ then implies that
$$
|K'|/|E(K')| + \delta(\d G) \le 1 + 2\epsilon\,.
$$
Since $G(K')$ is connected and $|K'| \to\infty$ when $\epsilon \to 0$,
it follows that $\beta(G) + \delta(\d G) \le 1$.

To prove that $\beta(G) + \delta(\d G) \ge 1$, let $\epsilon > 0$.
Let $K \subset V$ be
connected and simply connected (when regarded as a
union of closed faces of $\d G$ in the plane) such that
$|K|/|E(K)| \le \beta(G) + \epsilon$.
Let $\fK$ be the set of vertices in $\d G$
corresponding to the faces of $G(K)$.
Since
$|E^*(\fK)| \le |E(K)|$ and the number of faces of
the graph $G(K)$ is precisely $|\fK| + 1$,
we have
$$
|K|/|E(K)| + |\fK|/|E^*(\fK)|
\ge
|K|/|E(K)| + |\fK|/|E(K)|
= 
1 + 1/|E(K)|
\ge 1
$$
by Euler's formula applied to the graph $G(K)$.
(In case $\fK$ is empty, a comparable calculation shows that $|K|/|E(K)|
\ge 1$.)
In light of \rref e.deltasup/, it follows that
$$
\beta(G) + \delta(\d G) + \epsilon
\ge
\beta(G) + \epsilon + |\fK|/|E^*(\fK)|
\ge
|K|/|E(K)| + |\fK|/|E^*(\fK)|
\ge 1
\,.
$$
Since $\epsilon$ is arbitrary, the desired inequality follows.
\Qed

{}From \rref e.betaisoe/ and \rref e.deltaisoe/, we see that $\beta(G) >
\delta(G)$ when $G$ is regular and $\isoe'(G) > 0$. Thus, we obtain the
following inequality.

\procl c.beta_beta
If $G$ is a planar regular graph with nonamenable regular dual, then
$$
\beta(G) + \beta(\d G) > 1\,.
$$
\endprocl

The proof of \rref t.beta_delta/ appears not to give any idea of which
finite sets $K \subset V$ yield quotients $|\bde K|/|K|$ close to
$\isoe'(G)$. However,
Y.~Peres has deduced the following from a closer
examination of the proof. As in the proof of \rref t.beta_delta/,
we shall write $K'$ for the set of vertices incident to the faces
corresponding to $K$, for both $K \subset V$ and for $K \subset \d V$.
Likewise, $\widehat K$ denotes the faces inside the outermost cycle of
$E(K')$. According to the reasoning of the first paragraph of the proof of
\rref t.beta_delta/ and \rref e.upperbnd/, we have
\rlabel upperbnd2
{|(\widehat K)'|/|E((\widehat K)')| + |K|/|E^*(K)| 
\le
|(\widehat K)'|/|E((\widehat K)')| + |\widehat K|/|E^*(\widehat K)| 
\le 1 + 1/|E((\widehat K)')| \,.
}

\procl p.opt_sets
Let $G$ be a planar regular graph with regular dual $\d G$.
Let $K_0 \subset V$ be an arbitrary finite connected set
and recursively define $L_n := (\widehat K_n)' \subset \d
V$ and $K_{n+1} := (\widehat L_n)' \subset V$.
Then $|\bde K_n|/|K_n| \to \isoe'(G)$
and $|\bde L_n|/|L_n| \to \isoe'(\d G)$. 
\endprocl


\proof
The three amenable cases are trivial, so assume that $G$ is nonamenable.
Write 
$$
\kappa_n :=  |\bde K_n|/|K_n| - \isoe'(G)
$$
and 
$$
\lambda_n := |\bde L_n|/|L_n| - \isoe'(\d G)
\,.
$$
Also write $d := d_G$, $\d d := d_{\d G}$, $\iota := \isoe'(G)$, and
$\d \iota := \isoe'(\d G)$.
We may rewrite \rref e.upperbnd2/ as
$$
{2 \over \d d-|\bde L_n|/|L_n|} + {2 \over d+|\bde K_n|/|K_n|}
\le 1 + {1 \over |E(L_n)|}
\,,
$$
or, again, as
$$
{2 \over \d d-\d\iota-\lambda_n} + {2 \over d+\iota+\kappa_n}
\le 1 + {1 \over |E(L_n)|}
=
{2 \over \d d-\d\iota} + {2 \over d+\iota} + {1 \over |E(L_n)|}
\,,
$$
whence
$$
{2 \lambda_n \over (\d d-\d\iota)(\d d-\d\iota-\lambda_n)}
+ {2 \kappa_n \over (d+\iota)(d+\iota+\kappa_n)}
\le {1 \over |E(L_n)|}
\,.
$$
Therefore
\begin{align*}
2\lambda_n 
&\le
{(\d d-\d\iota)(\d d-\d\iota-\lambda_n) \over
(d+\iota)(d+\iota+\kappa_n)} (2\kappa_n)
+ {(\d d-\d\iota)(\d d-\d\iota-\lambda_n)
\over |E(L_n)|}
\cr&\le
\left({\d d-\d\iota \over d+\iota}\right)^2 2\kappa_n
+ {(\d d-\d\iota)^2 \over |E(L_n)|} \,.
\end{align*}
Similarly, we have
$$
2\kappa_{n+1}
\le
\left({d-\iota \over \d d+\d\iota}\right)^2 2\lambda_n + {(d-\iota)^2 \over
|E(K_{n+1})|}
\,.
$$
Putting these together, we obtain
$$
2\kappa_{n+1} \le a(2\kappa_n) + b_n
\,,
$$
where
$$
a :=
\left({(d-\iota)(\d d-\d\iota) \over (d+\iota)(\d d+\d\iota)}\right)^2
$$
and
$$
b_n :=
\left({(d-\iota)(\d d-\d\iota) \over \d d+\d\iota}\right)^2
{1 \over |E(L_n)|} + {(d-\iota)^2 \over |E(K_{n+1})|}
\,.
$$
Therefore
$$
2\kappa_n \le 2\kappa_0 a^{n-1} + \sum_{j=0}^{n-2} a^j b_{n-j}
\,.
$$
Since $a < 1$ and $b_n \to 0$, we obtain $\kappa_n \to 0$. Hence $\lambda_n
\to 0$ too.
\Qed

The following proposition relating the geometry of the graph to the
behavior of the free random-cluster measure uses the ideas of \cite{J, JS}.

\procl p.freebound
Let $G$ be a graph with degrees bounded by $d$. Write $b :=
\log\bigl(p/(1-p)\bigr)/\log q$ and $b^+ := \max\{b, 0\}$.
If $b < \beta(G)$ and 
$$
\log q > \frac{1+\log (d-1)}{\beta(G) - b^+}
\,,
$$
then $\FRC^G_{p,q}$-a.s., there is no infinite cluster.
\endprocl

\procl r.anchored
In this proposition, a better result is obtained by replacing
$\beta(G)$ by the corresponding quantity that results when $K$ is
required to contain a fixed point $o$ in the definition of $\beta(G)$.
The same proof applies. 
\endprocl

To prove \rref p.freebound/, we shall use the following bound analogous to
the well-known bound of Kesten \cite{Kes} on site-connected clusters.

\procl l.bondanimals
Let $G$ be a graph with degrees bounded by $d$. For
any fixed $o \in V(G)$, let $b_n$ be the number of connected subgraphs of
$G$ that contain $o$ and have exactly $n$ edges. Then $\limsup_{n \to\infty}
b_n^{1/n} < e(d-1)$.
\endprocl

\proof
Let $b_{n, \ell}$ denote the number of connected subgraphs $(V', E')$
of $G$ such that $o \in V'$, $|E'|= n$, and $|E^*(V')\setminus E'| = \ell$.
Note that for such a subgraph, 
\rlabel bndsize
{\ell \le d|V'| - 2n \le d(n+1)-2n = (d-2)n+d
\,.
}
Let $p := 1/(d-1)$ and consider Bernoulli($p$) bond percolation on $G$.
Writing the fact that 1 is at least the probability that the cluster of $o$
is finite and using \rref e.bndsize/, we obtain
$$
1 \ge \sum_{n, \ell} b_{n, \ell} p^n (1-p)^\ell
\ge
\sum_n b_n p^n (1-p)^{(d-2)n+d}
\,.
$$
Therefore $\limsup_{n \to\infty} b_n^{1/n} \le
1/[p(1-p)^{d-2}]$. Putting in the chosen value of $p$ gives the result.
\Qed

\proofof p.freebound
Because of \rref e.witness_genfree/, it suffices to prove the claim when $b
\ge 0$. So assume that $b \ge 0$.

Let $o \in V(G_i)$ for all $i$. If $\xi$ is a configuration, write $\xi(o)$
for the component of $o$ determined by $\xi$.
Suppose that $G' = \big(V(G'), E(G')\big)$ is a finite connected subgraph
of $G$ containing $o$.
If $\xi$ is a configuration in $G_i$ such that $\xi(o) =
G'$, then let $\xi'$ be obtained from $\xi$ by closing all edges in
$E(G')$. Then $\kk(\xi') = \kk(\xi) + |V(G')| - 1$, whence
$$
\RC^{G, i}_{p, q}[\xi]
=
\left({p \over 1-p}\right)^{|E(G')|} q^{-|V(G')|+1} \RC^{G, i}_{p, q}[\xi']
\,.
$$
Also,
if $\xi_1 \ne \xi_2$ are such that $\xi_1(o) = \xi_2(o) = G'$, then $\xi_1'
\ne \xi_2'$. Summing over all $\xi$ such that $\xi(o) = G'$ yields
$$
\RC^{G, i}_{p, q}[\xi(o) = G']
\le
\left({p \over 1-p}\right)^{|E(G')|} q^{-|V(G')|+1}
\,.
$$
Now the right-hand side equals $q^{b |E(G')|-|V(G')|+1}$ by the
definition of $b$.  Choose $b' \in \bigl(b, \beta(G)\bigr)$. Then
provided that $|V(G')|$ is sufficiently large, we have that the right-hand
side
is less than $q^{(b - b') |E(G')|+1}$. Therefore, for all large $N$, we
have, by \rref l.bondanimals/, that
$$
\RC^{G, i}_{p, q}\bigl[|E\bigl(\xi(o)\bigr)| \ge N\bigr]
\le
\sum_{n \ge N} \bigl(e(d-1)\bigr)^n q^{(b-b')n+1}
< \infty
$$
if $\log q > \bigl(1+\log (d-1)\bigr)/(b'-b)$. 
Therefore, for such $q$ this sum tends to 0 as $N
\to\infty$, so that given any $\epsilon > 0$, there is some $N_0$ such
that for all $N \ge N_0$ and for all $i$, we have $\RC^{G, i}_{p,
q}\bigl[|E\bigl(\xi(o)\bigr)| \ge N\bigr] < \epsilon$. Since the event
$\{|E\bigl(\xi(o)\bigr)| \ge N\}$
depends only on the edges within distance $N$ of $o$, it follows that
$\FRC^{G}_{p, q}\bigl[|E\bigl(\xi(o)\bigr)| \ge N\bigr] < \epsilon$. Since
$\epsilon > 0$ is
arbitrary, we find that $\FRC^{G}_{p, q}\bigl[|E\bigl(\xi(o)\bigr)| =
\infty\bigr] = 0$.
\Qed

\procl c.wiredbound
Let $G$ be a planar quasi-transitive graph whose dual has degrees
bounded by $\d d$.  Write $b := \log\bigl(p/(1-p)\bigr)/\log q$. If $b >
1 - \beta(\d G)$ and 
$$
\log q > \frac{1+\log (\d d-1)}{b\wedge1 - 1 + \beta(\d G)}
\,,
$$
then $\WRC^G_{p,q}$-a.s., there is a unique infinite cluster.
\endprocl

\proof
Let $p'$ be as in \rref e.p'/. Then $p'/(1-p') = q(1-p)/p$.
Let $b' := \log\bigl(p'/(1-p')\bigr)/\log q = 1 - b$. By our
hypothesis, $b' < \beta(\d G)$ and $\log q > \bigl(1+\log (\d
d-1)\bigr)/(\beta(\d G) - b^+)$, whence \rref p.freebound/ applied to $\d G$
shows that there is no infinite cluster $\FRC^{\d G}_{p', q}$-a.s.
Hence the conclusion follows from \rref p.FWcoupled/.
\Qed

Putting all this together, we arrive at our main result in this section:

\procl t.separated
If $G$ is a planar regular nonamenable graph of degree $d$ with regular
dual of degree $\d d$, then $\pcf(q) > \puw(q)$ for 
$$
q > {\bigl(2+\log\bigl((d-1) (\d d-1)\bigr)\bigr)(d \d d - d - \d
d) \over \sqrt{(d-2)(\d d-2)(d \d d - 2d - 2\d d)}}
\,.
$$
\endprocl

\proof
Let 
$$
b_0 := {\bigl(1+\log (d-1)\bigr)\bigl(1-\beta(\d G)\bigr) + \bigl(1+\log (\d
d-1)\bigr)\beta(G) \over
\bigl(1+\log (d-1)\bigr) + \bigl(1+\log (\d d-1)\bigr)}
\,.
$$
Then $0 < 1 - \beta(\d G) < b_0 < \beta(G) < 1$ because of \rref
c.beta_beta/. Furthermore, we have
$$
{1 + \log (d-1) \over \beta(G)-b_0^+}
=
{1 + \log (\d d-1) \over b_0\wedge1-1+\beta(\d G)}
=
{\bigl(2+\log\bigl((d-1) (\d d-1)\bigr)\bigr)(d \d d - d - \d d) \over
\sqrt{(d-2)(\d
d-2)(d \d d - 2d - 2\d d)}}
$$
since 
$$
\beta(G) =
{d(\d d-2) + \sqrt{(d-2)(\d d-2)(d\d d - 2d - 2\d d)} \over 2(d \d d -
d - \d d)}
\,,
$$
which follows from \rref t.isocalc/ and some calculation.
Thus, there is a positive-length
interval of $p$ for which $b$ in \rref p.freebound/ is
close enough to $b_0$ that the hypotheses of both
\rref p.freebound/ and \rref c.wiredbound/ are satisfied. This gives the
result.
\Qed

\procl r.when_equal
Even when $\pcf=\puw$, there is no infinite cluster $\FRC_{p,q}$-a.s.\ and
a unique infinite cluster $\WRC_{p,q}$-a.s.\ for $p = \pcf=\puw$, in view
of \rref t.death/ and \rref c.uniq_p_u/.
\endprocl

For comparison, we present the following bounds on when the opposite
inequality $\pcf(q) < \puw(q)$ holds.

\procl p.small_q
If $G$ is a planar regular nonamenable graph of degree $d$ with regular
dual of degree $\d d$, then $\pcf(q) < \puw(q)$ for $q < d \d d - 2d - 2\d
d$.
\endprocl

For example, if $d = \d d = 5$, then $\pcf<\puw$ if $q<5$, while
$\pcf>\puw$ if $q > (2+4\log2)\sqrt5 = 10.67^+$.

\proof
We claim that $\pcf(G, q) < \puw(G, q)$ if 
\rlabel{sep}
{
h\big(\pc(G)\big) h\big(\pc(\d G)\big) < 1/q
\,,
}
where $h$ is defined as in \rref e.hdef/.
Indeed, 
\rlabel{comp1}
{
h\big(\pcf(G,q)\big) \le q h\big(\pc(G)\big)
}
by \rref e.witness_genfree/
(applied with $p_1 = p$, $q_1 = 1$, $q_2 = q$, and $p_2$ chosen so that
$h(p) = h(p_2)/q$) and
\rlabel{comp2}
{
h\big(\puw(G,q)\big) = q/h\big(\pcf(\d G,q)\big) \ge 1/h\big(\pc(\d G)\big)
}
by \rref c.pcpu/ and \rref e.comp1/ (applied to $\d G$). Therefore, if
\rref e.sep/ holds, then \rref e.comp1/ and \rref e.comp2/ give
$$
h\big(\pcf(G,q)\big) \le q h\big(\pc(G)\big)
< 1/h\big(\pc(\d G)\big)
\le h\big(\puw(G,q)\big)
\,,
$$
which implies that $\pcf(G, q) < \puw(G, q)$.

Now $\pc(G) \le 1/(1+\isoe(G))$ by \cite{BS}. Therefore \rref e.sep/ is
implied by $\isoe(G) \isoe(\d G) > q$. This is the same as the claimed
range of $q$, as \rref t.isocalc/ and some calculation shows.
\Qed

\procl r.ising_extreme
Note that when $\pcf(G) < \puw(G)$, we also have $\pcw(G) < \puf(G)$, and
thus we have an interval of $p$ for which $\FRC_p \ne \WRC_p$ (see \rref
p.f=w/ above). For $q=2$, the condition in \rref p.small_q/ shows that this
holds for $\d d = 3$ when $d > 8$, for $\d d = 4$ when $d > 5$, and for all
$\d d \ge 5$ when $d \ge 5$. This generalizes the main result in \cite{Wu}. 
\endprocl

We end the section with some open questions.

\procl q.pupcphase
Let $G$ be a nonamenable quasi-transitive graph.
Is the set of $q$ for which $\pcf(G, q) > \puw(G, q)$ an interval?
If $G$ has only one end, is the set of such $q$ nonempty?
\endprocl

\procl q.distinct
If $G$ is nonamenable and quasi-transitive, can there be there any $q
> 1$ such that $\pcw(G,q) = \puf(G,q)$ (so that all four critical values
coincide)?
\endprocl

\section{Robust phase transition} \label{sect:robust}

A classical question about the Potts model is when it exhibits
a phase transition for given $q$, $\beta$ and $G$.
We say that a {\bf phase transition} occurs when there is more than
one Gibbs measure for the $q$-state Potts model on $G$ with 
inverse temperature $\beta$.
As noted above in Section \ref{sect:critical}, this happens iff
$\WPt^G_{q,\beta,r_1} \neq \WPt^G_{q,\beta,r_2}$ for
$r_1 \neq r_2$.
This is equivalent to the statement $\WPt^G_{q,\beta,r}(\omega(o)=r)>1/r$
for any fixed vertex $o \in V$ and also
to $\WRC^G_{p,q}(o \leftrightarrow \infty) > 0$, where
$p := 1-e^{-2\beta}$ and $\{o \leftrightarrow \infty\}$ is the event
that $o$ is contained in an infinite open cluster.
Let $\{o \leftrightarrow \partial V_i\}$ be the event that
$o$ is connected to $\partial V_i$ by a path of open edges in $G_i$
(where $G_i$, $V_i$ and $E_i$ are as in Section \ref{sect:background}).
It follows that phase transition in the $q$-state Potts model
with the given parameters is equivalent to
$\inf_i \WRC^{G,i}_{p,q}(o \leftrightarrow \partial V_i)>0$.
Hence there exists a critical value $\betac$ 
(given by $\betac= -\frac{1}{2}\log(1-\pc^\w)$)
such that we have phase transition for $\beta>\betac$, but not for
$\beta<\betac$.

Pemantle and Steif \cite{PS} introduced the stronger concept
of robust phase transition. Although they considered mostly the 
Heisenberg model, they also have some results concerning robust 
phase transition in the Potts model. 
In order to define this notion, we need to generalize the Potts model slightly:
When defining $\Pt^{G_i}_{q,\beta}$, let us allow different 
interaction along different edges, i.e., replace
$\beta$ with ${\mathbf B} := \{\beta_e\}_{e \in E_i}$.
It is then straightforward to modify the measures $\Pt^{G_i}_{q,{\mathbf
B}}$
to measures $\FPt^{G,i}_{q,{\mathbf B}}$ and $\WPt^{G,i}_{q,{\mathbf B}}$
in the same way as we did in Section \ref{sect:background}
for the case ${\mathbf B} \equiv \beta$.
Now for $\epsilon>0$, define $\B^i_{\epsilon}$ so that
$B^i_{\epsilon}(e)$ equals $\epsilon$ for all edges with one endpoint
in $\partial V_i$ and equals $\beta$ for all other edges.
We say that the $q$-state Potts model with inverse temperature $\beta$
exhibits a {\bf robust phase transition} if
$\inf_i\phantom{\Big|}
\WPt^{G,i}_{q,\B^i_{\epsilon},r}\bigl(\omega(o)=r\bigr)>1/r$
for {\em all} $\epsilon>0$.

By making a corresponding extension of the random-cluster model,
i.e., by replacing $p$ with $\p := \{p_e\}_{e \in E_i}$,
we see that this is equivalent to
\[\inf_i \WRC^{G,i}_{\p^i_s,q}(o \leftrightarrow \partial V_i) > 0\]
for all $s>0$, where $\p^i_s$ equals $p:=1-e^{-2\beta}$ for edges that do
not have an endpoint in $\partial V_i$ and equals $s$ for those that do.
Pemantle and Steif show that when $G$ is a tree and $q \geq 3$, then there
is sometimes a phase transition but not a robust phase transition.
In particular, they show the following.  
On trees, the critical value 
$\betac^{{\rm robust}}$ for robust phase transition is a strictly decreasing 
function of the branching number of the tree, but $\betac$ is not.
Consequently, if $T_1$ and $T_2$ are two trees with $\br(T_1)<\br(T_2)$ but where
the $q$-state Potts model exhibits a phase transition on $T_1$ but not on 
$T_2$, then this happens for an interval of inverse temperatures on $T_1$.
For instance, if $G$ is taken to be the binary tree
and $q \geq 3$, then there exists an $\epsilon>0$ such that
there is a phase transition but not a robust phase transition
when $e^{2\beta}-1 \in (q-\epsilon,q)$. See Theorems 1.13 
and 1.14 of \cite{PS}.
We show:

\procl t.robust
Let $G$ be an infinite regular nonamenable graph with degree $d$.
Then there exists $q_0<\infty$ such that for $q \geq q_0$ and
$e^{2\beta}-1 \in [q^{2/(d+\iota'_{E}(G)/2)},q^{2/(d-\iota'_{E}(G)/2)}]$,
the $q$-state Potts model on $G$ with inverse temperature $\beta$
exhibits a phase transition but not a robust phase transition.
\endprocl

\proof
Fix $i$. From the definition of $\WRC^{G,i}_{\p^i_s,q}$ it is clear that
we may, instead of regarding the vertices outside 
$V_i$ as connected, regard them
as contracted to a single vertex $v_0$.
Now define a new graph $H_i$ by adding an edge between $v_0$ and
all vertices of $V_i$, i.e.,
let $V(H_i) := V_i \cup \{v_0\}$
and $E(H_i) := E(V_i) \cup
E_0$ where $E_0 := \{\ed{v_0,v} \st v \in V_i\}$.
Let $\tilde{\p}$ be $s$ for all edges of $E_0$
and $p$ for all other edges.
Then by conditioning on
$\RC^{H_i}_{\tilde{\p},q}(E_0)$ and using Holley's inequality, we see that
\rlabel compare-to-old
{\RC^{H_i}_{\tilde{\p},q}\bigl(E^*(V_i)\bigr) \geqd \WRC^{G,i}_{\p^i_s,q}.}
Now the proof of \cite{J}, Theorem 1.2(a),
shows precisely that for $q$ large enough, $s$ 
small enough, and $\beta$ in the range indicated,
$\inf_i \RC^{H_i}_{\tilde{\p},q}(A_i) = 0$,
where $A_i$ is the event that there is a path of open
edges connecting $o$ to $\partial V_i$ without passing through $v_0$.
[The proof is essentially the same as that of \rref p.freebound/,
but here one has to take the edges of $E_0$ into account.
This is done by a simple application of Markov's inequality.]
By \rref e.compare-to-old/, it follows that
there is no robust phase transition in the corresponding
Potts model.
On the other hand, \cite{J}, Theorem 4.4(a), says that for $q$ large enough,
a phase transition occurs for $\beta$ in the range indicated.
\Qed

Let us now briefly consider the case when $G$ is instead an amenable
quasi-transitive graph.
Then, as noted above, $\FRC_{p,q} = \WRC_{p,q}$ for all but at most
countably
many values of $p$.
Therefore, if $\WRC_{p,q}(o \leftrightarrow \infty)>0$, then
$\FRC_{p',q}(o \leftrightarrow \infty) >0$ for all $p'>p$.
This strongly suggests that 
$\inf_i \FRC^{G,i}_{p',q}(o \leftrightarrow \partial V_i)>0$,
and since $\FRC^{G,i}_{p,q}$ is $\WRC^{G,i}_{\p^i_s,q}$
with $s=0$, this would imply that the critical temperatures for
phase transition and robust phase transition coincide.
However, as $\FRC^{G,i}_{p,q}$ is increasing in $i$, it is
possible to imagine a scenario where $o$ is not connected
to $\partial V_i$ for any $i$ but still connected to infinity
in the limit.
This does not happen for the $\WRC$ measures on any graph (\cite{ACCN},
proof of Theorem 2.3(c)), but is unknown for $\FRC$.

\procl q.amenable_robust
If $G$ is amenable and quasi-transitive, then is the critical
inverse temperature for robust phase transition $\betac$?
\endprocl

This question does not address the case
$\beta=\betac$.
For $G := {\mathbf Z}^d$, $d \geq 2$ and $q$ large, it is well known
that the Potts model at criticality ($\beta=\betac$) exhibits phase
transition, whereas it was recently shown by 
van Enter \cite{vE} that (still at criticality) there is no 
robust phase transition.

\medbreak
\noindent {\bf Acknowledgement.}\enspace
We are grateful to
Yuval Peres for permission to include \rref p.opt_sets/ and its proof,
to Roberto Schonmann for a careful reading and several useful
comments, and to Oded Schramm for discussions.


\begin{thebibliography}{99}

{\small

\bibitem{ACCN} Aizenman, M., Chayes, J.T., Chayes, L., and Newman, C.M.
(1988) Discontinuity of the magnetization in one-dimensional
$1/|x-y|^2$ Ising and Potts models, {\sl J. Statist. Phys.} {\bf 50},
1--40.

\bibitem{BB} Babson, E. and Benjamini, I. (1999) Cut sets and normed
cohomology with applications to percolation, {\sl Proc. Amer. Math.
Soc.} {\bf 127}, {589--597}.

\bibitem{BauesPeyerimhoff} Baues, O. and Peyerimhoff, N. (2001) Curvature
and geometry of tessellating plane graphs, {\sl Discrete Comput. Geom.}
{\bf 25}, 141--159.

\bibitem{BLPS} Benjamini, I., Lyons, R., Peres, Y., and Schramm, O. (1999)
Group-invariant percolation on graphs, {\sl Geom. Funct. Analysis} 
{\bf 9}, 29--66.

\bibitem{BLPScrit} Benjamini, I., Lyons, R., Peres, Y., and Schramm, O.
(1999).  Critical percolation on any nonamenable group has no infinite
clusters.  {\sl Ann. Probab.} {\bf 27}, 1347--1356.

\bibitem{BS} Benjamini, I. and Schramm, O. (1996) Percolation beyond
${\bf Z}^d$, many questions and a few answers, {\sl Electr. Commun. Probab.}
{\bf 1}, 71--82. 

\bibitem{BS:hp} Benjamini, I. and Schramm, O. (2001) Percolation in
the hyperbolic plane, {\sl J. Amer. Math. Soc.} {\bf 14}, 487--507.  

\bibitem{BRZ} Bleher, P.M., Ruiz, J. and Zagrebnov, V.A. (1995) On the
purity
of the limiting Gibbs state for the Ising model on the Bethe lattice,
{\sl J. Statist. Phys.} {\bf 79}, 473--482. 

\bibitem{BBCK} Biskup, M., Borgs, C., Chayes, J.T., and Koteck\'y, R.
(2000) Gibbs states of graphical representations in the {P}otts model with
external fields,  {\sl J. Math. Phys.} {\bf 41}, 1170--1210.

\bibitem{BC:covariance} Borgs, C. and Chayes, J.T. (1996) The covariance
matrix of the {P}otts model: a random cluster analysis, {\sl J. Statist.
Phys.} {\bf 82}, 1235--1297.

\bibitem{BK} Burton, R.M. and Keane, M.S. (1989) Density and uniqueness in
percolation, {\sl Commun. Math. Phys.} {\bf 121}, 501--505.

\bibitem{CCS}  Chayes, J.T., Chayes, L., and Schonmann, R.H. (1987)
Exponential decay of connectivities in the two-dimensional Ising model,
{\sl J. Statist. Phys.} {\bf 49}, 433--445.

\bibitem{CCST} Chayes, J.T., Chayes, L., Sethna, J.P., and Thouless, D.J.
(1986).  A mean field spin glass with short-range interactions,  {\sl Comm.
Math. Phys.} {\bf 106}, 41--89.

\bibitem{CK} Chaboud, T. and Kenyon, C. (1986)
Planar Cayley graphs with regular dual,
{\sl Internat. J. Algebra Comput.} {\bf 6}, no. 5, 553--561.

\bibitem{ES} Edwards, R.G. and Sokal, A.D. (1988) Generalization of the
Fortuin-\-Kasteleyn-\-Swendsen-\-Wang representation and Monte Carlo
algorithm, {\sl Phys. Rev.} D {\bf 38}, 2009--2012.

\bibitem{vE} van Enter, A. (2000) A remark on the notion of robust phase
transitions, {\sl J. Statist. Phys.} {\bf 98}, 1409--1417. 

\bibitem{EKPS} Evans, W., Kenyon, C., Peres, Y., and Schulman, L. (2000)
Broadcasting on trees and the Ising model, 
{\sl Ann. Appl. Probab.} {\bf 10}, 410--433. 

\bibitem{FP} Floyd, W.J. and Plotnick, S.P. (1987)
Growth functions on Fuchsian groups and the Euler characteristic, {\sl
Invent. Math.} {\bf 88}, 1--29.

\bibitem{FK} Fortuin, C.M. and Kasteleyn, P.W. (1972)
On the random-cluster model. I. Introduction and relation to other models,
{\sl Physica} {\bf 57}, 536--564.

\bibitem{GKN} Gandolfi, A., Keane, M.S. and Newman, C.M. (1992)
Uniqueness of the infinite component in a random graph with
applications to percolation and spin glasses,
{\sl Probab. Theory Related Fields} {\bf 92}, 511--527.

\bibitem{Geor} Georgii, H.-O. (1988)
{\sl Gibbs Measures and Phase Transitions.}
Walter de Gruyter \& Co., Berlin-New York.

\bibitem{GHM} Georgii, H.-O., H\"aggstr\"om, O. and Maes, C. (2001)
The random geometry of equilibrium phases, {\sl Phase Transitions and
Critical Phenomena} (C. Domb and J.L. Lebowitz, eds), pp. 1--142, 
Academic Press, London.

\bibitem{GH} Grigorchuk, R. and de la Harpe, P. (1997)
On problems related to growth, entropy, and spectrum in group theory,
{\sl J. Dynam. Control Systems} {\bf 3}, 51--89. 

\bibitem{Gr} Grimmett, G.R. (1995) The stochastic random-cluster process,
and
the uniqueness of random-cluster measures, {\sl Ann. Probab.} {\bf 23}, 
1461--1510. 

\bibitem{Gr:compare} Grimmett, G.R.  (1995)
Comparison and disjoint-occurrence inequalities for random-cluster models,
{\sl J. Statist. Phys.} {\bf 78}, 1311--1324.

\bibitem{GN} Grimmett, G.R. and Newman, C.M. (1990) Percolation in
$\infty+1$
dimensions, {\sl Disorder in Physical Systems} (G.R. Grimmett and D.J.A. 
Welsh, eds), 167--190, Clarendon Press, Oxford. 

\bibitem{H2} H\"aggstr\"om, O. (1996) The random-cluster model on a
homogeneous tree, {\sl Probab. Th. Rel. Fields} {\bf 104}, 231--253.

\bibitem{H3} H\"aggstr\"om, O. (1997) Infinite clusters in dependent 
automorphism invariant percolation on trees, {\sl Ann. Probab.} {\bf 25}, 
1423--1436. 

\bibitem{H} H\"aggstr\"om, O. (1998) Random-cluster representations in
the study of phase transitions, {\sl Markov Proc. Rel. Fields} {\bf 4},
275--321.

\bibitem{HJL} H\"aggstr\"om, O., Jonasson, J. and Lyons, R. (2001) Coupling
and Bernoullicity in random-cluster and Potts models, {\sl preprint}.

\bibitem{J} Jonasson, J. (1999)
The random cluster model on a general graph and a phase
transition characterization of nonamenability, {\sl Stoch. Proc. Appl.}
{\bf 79}, 335--354.

\bibitem{JS} Jonasson, J. and Steif, J.E. (1999)
Amenability and phase transition in the Ising model,
{\sl J. Theor. Probab.} {\bf 12}, 549--559. 

\bibitem{Kes} Kesten, H. (1982) {\sl Percolation Theory for
Mathematicians.} Birkh\"auser, Boston, Mass.

\bibitem{L} Lyons, R. (2000) Phase transitions on nonamenable graphs, 
{\sl J. Math. Phys.} {\bf 41}, 1099--1126. Version with correct proof of
tail triviality of random-cluster measures available at
\texttt{arXiv:math.PR/9908177}.

\bibitem{LP:book} 
{Lyons, R. with Peres, Y.} (2001)
{\sl Probability on Trees and Networks}.
Cambridge University Press.
In preparation. Current
  version available at \hfil\break
  {\tt http://php.indiana.edu/\string~rdlyons/}.

\bibitem{LS} Lyons, R. and Schramm, O. (1999) Indistinguishability of 
percolation clusters, {\sl Ann. Probab.} {\bf 27}, 1809--1836.

\bibitem{LSS} Liggett, T.M., Schonmann, R.H. and Stacey, A.M. (1997)
Domination by product measures, {\sl Ann.  Probab.} {\bf 25}, 71--95.

\bibitem{NS} Newman, C.M. and Schulman, L.S. (1981) Infinite clusters in
percolation models, {\sl J. Statist. Phys.} {\bf 26}, 613--628. 

\bibitem{PS} Pemantle, R. and Steif, J.E. (1999) Robust phase transition for
Heisenberg and other models on general trees, {\sl Ann. Probab.}
{\bf 27}, 876--912.

\bibitem{PW} Propp, J.G. and Wilson, D.B. (1996)
Exact sampling with coupled Markov chains and applications to statistical
mechanics, {\sl Random Structures Algorithms} {\bf 9}, 223--252.

\bibitem{Schon:mult} Schonmann, R.H. (2000) Multiplicity of phase
transitions and mean-field criticality on highly non-amenable graphs, {\sl
Comm. Math. Phys.}, to appear.

\bibitem{SeriesSinai}
Series, C.M. and Sina{\u\i}, Y.G. (1990)
Ising models on the {L}obachevsky plane,
{\sl Comm. Math. Phys.} {\bf 128}, 63--76.

\bibitem{SW} Swendsen, R.H. and Wang, J.-S. (1987)
Nonuniversal critical dynamics in Monte Carlo simulations,
{\sl Phys. Rev. Lett.} {\bf 58}, 86--88.

\bibitem{Tr} Trofimov, V.I. (1985) Automorphism groups of graphs as
topological groups, {\sl Math.\ Notes} {\bf 38}, 717--720. 

\bibitem{Wel} Welsh, D.J.A. (1993) Percolation in the random cluster process
and the $Q$-state Potts model, {\sl J. Phys. A} {\bf 26}, 2471--2483.

\bibitem{Wu} Wu, C.C. (2000) Ising models on hyperbolic graphs II, {\sl
J. Stat. Phys.} {\bf 100}, 893--904.

}
\end{thebibliography}
\end{document}